\documentclass[11pt]{amsart}
\usepackage[dvips]{epsfig}
\usepackage{graphics}
\usepackage{latexsym}
\usepackage{verbatim}
\usepackage{amsmath}
\usepackage{amsthm}
\usepackage{amssymb}
\newtheorem{theorem}{Theorem}[section]
\newtheorem{lemma}[theorem]{Lemma}
\newtheorem{corollary}[theorem]{Corollary}
\newtheorem{proposition}[theorem]{Proposition}
\newtheorem{definition}[theorem]{Definition}

\def\Remark{\medskip\noindent{\bf Remark: }}
\def\Remarks{\medskip\noindent{\bf Remarks: }}

\newcommand{\ens}[1]{\mathbb{#1}}

\newcommand{\N}{\mathbb{N}}

\newcommand{\R}{\mathbb{R}}

\def\cal{\mathcal}
\def\cst{{\rm cst}} 

\def\derpar#1#2{\frac{\partial#1}{\partial#2}}
\def\var{\varepsilon}
\def\signcm{\bigskip\bigskip\hspace{80mm}
\vbox{{\sc C. Mouhot\par\vspace{3mm}
UMPA, \'ENS Lyon\par
46 all\'ee d'Italie\par
69364 Lyon Cedex 07\par
FRANCE\par\vspace{3mm}
e-mail:} cmouhot@umpa.ens-lyon.fr }}

\begin{document}

\title[Quantitative lower bounds for the full Boltzmann equation]
{Quantitative lower bounds for the full Boltzmann equation,  
Part I: Periodic boundary conditions}

\author{Cl\'ement Mouhot}

\hyphenation{bounda-ry rea-so-na-ble be-ha-vior pro-per-ties
cha-rac-te-ris-tic}

\begin{abstract} We prove the appearance of an explicit lower bound on the
solution to the full Boltzmann equation in the torus for a broad family of  
collision kernels including in particular long-range interaction models, 
under the assumption of some uniform bounds on some hydrodynamic quantities. 
This lower bound is independent of time and space. When the collision kernel 
satisfies Grad's cutoff assumption, the lower bound is a global Maxwellian and 
its asymptotic behavior in velocity is optimal, whereas for non-cutoff collision kernels 
the lower bound we obtain decreases exponentially but faster than the Maxwellian. 
Our results cover solutions constructed in a spatially homogeneous setting, 
as well as small-time or close-to-equilibrium solutions to the full Boltzmann 
equation in the torus. 
The constants are explicit and depend on the {\it a priori} bounds on the solution.  
\end{abstract}

\maketitle

\tableofcontents


\section{Introduction}
\setcounter{equation}{0}


This paper is devoted to the study of qualitative properties of the solutions 
to the full Boltzmann equation in the torus 
for a broad family of collision kernels. In this work we shall quantify the 
positivity of the solution by proving the ``immediate'' 
appearance of a stationary lower bound, uniform in space. Before we explain our 
results and methods in more details let us introduce the problem in a precise way. 

\subsection{Motivation}

The Boltzmann equation decribes the behavior of
a dilute gas; when we assume periodic boundary conditions in space, it reads
 \begin{equation*}\label{eq:base}
 \derpar{f}{t}  + v \cdot \nabla _x f = Q(f,f), \qquad x \in \ens{T}^N, \quad  v \in \R^N, \quad t \in [0,T),
 \end{equation*}
where $T \in (0, +\infty]$, $\ens{T}^N$ is the $N$-dimensional torus, 
the unknown $f=f(t,x,v)$ is a time-dependent
probability density on $\ens{T}^N _x \times \R^N _v$ ($N\geq2$), and $Q$ is the quadratic
Boltzmann collision operator. It is local in $t$ and $x$ and 
we define it by the bilinear form
 \begin{equation*}\label{eq:collop}
 Q(g,f) = \int _{\R^N} \, dv_*
 \int _{\ens{S}^{N-1}} \, d\sigma B(|v-v_*|, \cos \theta)\, (g'_* f' - g_* f).
 \end{equation*}
Here we have used the shorthands $f'=f(v')$, $g_*=g(v_*)$ and
$g'_*=g(v'_*)$, where
 \begin{equation*}\label{eq:rel:vit}
 v' = \frac{v+v_*}{2} + \frac{|v-v_*|}{2} \sigma, \qquad
 v'_* = \frac{v+v^*}{2} - \frac{|v-v_*|}{2} \sigma
 \end{equation*}
stand for the pre-collisional velocities of particles which after
collision have velocities $v$ and $v_*$. Moreover $\theta\in
[0,\pi]$ is the deviation angle between $v'-v'_*$ and $v-v_*$, and
$B$ is the Boltzmann collision kernel (related to the
cross-section $\Sigma(v-v_*,\sigma)$ by the formula $B=\Sigma
|v-v_*|$), determined by physics. On physical grounds, it is
assumed that $B \geq 0$ and that $B$ is a function of $|v-v_*|$ and
$\cos\theta= \left\langle \frac{v-v_*}{|v-v_*|}, \sigma \right\rangle$.


In this paper we shall be concerned with the case when $B$ takes the 
following product form 
 \begin{equation}\label{eq:B}
 B (v-v_*, \sigma) = \Phi(|v-v_*|) \, b(\cos \theta)
 \end{equation}
where $\Phi$ satisfies either the assumption 
 \begin{equation}\label{eq:Phi}
 \forall \, z \in \R, \ \ \ c_\Phi \; |z|^\gamma \le \Phi (z) \le C_\Phi \; |z|^\gamma
 \end{equation}
or the mollified assumption 
 \begin{equation}\label{eq:Phi:mol}
 \left\{
 \begin{array}{l}
 \forall \, |z| \ge 1, \ \ \  c_\Phi \; |z|^\gamma \le \Phi (z) \le C_\Phi \; |z|^\gamma \\
 \forall \, |z| \le 1, \ \ \  c_\Phi \le \Phi(z) \le C_\Phi 
 \end{array}
 \right. 
 \end{equation}
with $c_\Phi,C_\Phi >0$ and $\gamma \in (-N,1]$, and $b$ is a continuous function on $\theta \in (0,\pi]$, 
strictly positive on $\theta \in (0,\pi)$, such that 
 \begin{equation}\label{eq:b}
 b(\cos \theta) \sin^{N-2} \theta \ \sim_{\theta \to 0^+} \ b_0 \; \theta^{-1-\nu}
 \end{equation}
for some $b_0 > 0$ and $\nu \in (-\infty,2)$. 
The assumption~\eqref{eq:B} that $B$ takes a tensorial form is made for the sake of convenience, since it is 
a well-accepted hypothesis which covers important physical cases. Most probably 
it could be relaxed to non-tensorial collision kernels with the same kind of controls, up to some 
technical refinements. The assumption that $b$ is strictly positive on $\theta \in (0,\pi)$ 
is a technical requirement for Lemma~\ref{lem:up} and Lemma~\ref{lem:spread} in Section~\ref{sec:tools} 
and it could be relaxed to the requirement  
that $b$ is strictly positive near $\theta \sim \pi/2$.

Following the usual taxonomy we shall denote by ``hard potential'' collision kernels 
the case when $\gamma > 0$, ``Maxwellian'' collision kernels the case 
when $\gamma =0$, and ``soft potential'' collision kernels
the case when $\gamma <0$. When $\Phi$ satisfies assumption~\eqref{eq:Phi:mol} and 
not~\eqref{eq:Phi} we shall speak of ``mollified soft potentials'' collision kernels or 
``mollified hard potentials'' collision kernels. 
In the case when $\nu < 0$, the angular collision 
kernel is locally integrable, an assumption which is usually referred to as
{\em Grad's cutoff assumption} (see~\cite{Grad:58}). Thus for $\nu < 0$ we shall speak of 
``cutoff'' collision kernels, and for $\nu \ge 0$ we shall speak of ``non-cutoff'' collision kernels. 
  
The main cases of application are ``hard spheres'' interaction ($b$ constant and $\Phi(z)=|z|$ which 
corresponds to the case $\gamma =1$ and $\nu = -1$), 
and interactions deriving from a $1/r^s$ force ($s>2$), where $r$ is the distance between particles, which 
corresponds to $\gamma = (s - 5)/(s-1)$ and $\nu=2/(s-1)$ in dimension $3$. 

In the case when the solution $f(t,v)$ does not depend on the space variable $x$, we shall speak of 
spatially homogeneous solution. 
\medskip

The attempts to quantify the strict positivity of the solution to the Boltzmann equation 
are as old as the mathematical theory of the Boltzmann equation, since Carleman himself established 
such estimates in his pioneering paper~\cite{Carl:EB:32}. He showed,  
for hard potentials with cutoff in dimension $3$, 
that the spatially homogeneous solutions radially symmetric in $v$ 
he had constructed in $L^\infty _6 (\R^N _v)$~\footnote{see Subsection~\ref{subsec:not} for the notations}
(the very first result in the Cauchy theory) satisfy a lower bound of the form 
 \[ \forall \, t \ge t_0 >0, \ \forall \, v \in \R^3 \ \ \ f(t,v) \ge C_1 \, e^{-C_2 \, |v|^{2+\var}}, \] 
for any fixed $t_0 >0$ and $\var > 0$. The constants $C_1, C_2 >0$ are uniform as $t \to +\infty$ and 
depend on $t_0$, $\var$ and some estimates on the solution. The proof was based on a ``spreading property'' 
of the collision operator and the assumption that the initial datum is uniformly bounded from below 
by a positive quantity on a ball centered at the origin (in fact the weaker 
assumption of a lower bound on an annulus is sufficient, see~\cite{Carl:EB:32}). 

This result remained unchallenged until the paper from Pulvirenti and Wennberg~\cite{PW:lb:97}. 
They proved, for hard potentials with cutoff in dimension $3$, that the spatially homogeneous solutions  
in $L^1 _2 (\R^N _v)$\footnote{id.} with bounded entropy 
(see~\cite{Arke:I+II:72} and~\cite{MiscWenn:exun:99}) satisfy 
a lower bound of the form 
 \[ \forall \, t \ge t_0 >0, \ \forall \,  v \in \R^3 \ \ \ f(t,v) \ge C_1 \, e^{-C_2 \, |v|^2}, \] 
for any fixed $t_0 >0$. Again $C_1, C_2 >0$ are uniform as $t \to +\infty$ and depend on $t_0$  
and some estimates on the solution. Their proof was also based on the spreading property of 
the collision operator but the optimal decay of the lower bound was obtained by some  
improvements of the computations. Moreover they made a clever use of the iterated gain part of the collision 
kernel in order to establish the immediate appearance of a positive minoration of the solution on a ball, 
thus getting rid of the assumption of Carleman on the initial datum. 
This paper is the starting point of our study. 

Two other methods should be mentioned.

On one hand, Fournier~\cite{Four:Kac:00,Four:Kac:02} 
established by delicate probabilistic techniques that the spatially homogeneous solutions to 
the Kac equation without cutoff satisfy $f(t,v) \in C^\infty ((0,+\infty) \times \R^N)$ and  
 \[ \forall \, t > 0, \ \forall \, v \in \R^N, \ \ \ f(t,v) > 0. \]
He proved the same kind of result in~\cite{Four:Bol:01} for the spatially homogeneous 
Boltzmann equation without cutoff in dimension $2$ under technical restrictions. 

On the other hand there is a work in progress by Villani in order 
to establish lower bounds on the solution to the Boltzmann equation using suitable maximum principles. 
The most important feature of this new method is that it applies to the non-cutoff case. For more explanations 
we refer to~\cite[Chapter~2, Section~6]{Vill:hand}. 
This method has been able to recover more simply the results 
by Fournier, and a quantitative lower bound is in progress. We also refer to the 
work~\cite{PaVi:np} which proves with the same tools the propagation of upper Maxwellian bounds for the 
spatially homogeneous solutions to the Boltzmann equation for hard spheres. 

Finally we note that in the case of the spatially homogeneous Landau equation with Maxwellian or 
hard potentials interactions, one can prove a theorem similar to that of Pulvirenti and 
Wennberg by means of the standard maximum principle for 
parabolic equations, see~\cite{DV:LI:00}. Actually the result stated in this paper is not uniform as 
$t \to +\infty$, but it can be made uniform by the same argument as in the proof of 
Theorem~\ref{theo:main} in Section~\ref{sec:proof} below. 
\medskip

The study of lower bounds is of interest in itself, in order to understand 
the qualitative behaviour of solutions to the Boltzmann equation. 
Moreover recently this interest has been renewed by 
the emergence of a new quantitave method in the study of convergence to 
equilibrium, the so-called ``entropy-entropy production'' method 
(see~\cite{ToVi:EEP:99,ToscVill:soft:00,DV:LII:00,Vill:EEP:03}). 
This method requires indeed a control from below on the solution by a function decreasing 
at most exponentially, and uniform in time. 
It has been applied lately to some inhomogeneous kinetic equations: see~\cite{DV:eqFP:01} for 
the Fokker-Planck equation and~\cite{DV:eqEB:03} for the full Boltzmann equation. 
For instance the main result in~\cite{DV:eqEB:03} asserts that any solution of the Boltzmann equation 
satisfying uniform estimates of smoothness and fast decay at large velocities, combined with a lower 
bound like 
 \begin{equation} \label{eq:hypEEP}
 \forall \, t \ge t_0 >0, \ \forall \, x \in \ens{T}^N, \ \forall \,  v \in \R^N, \ \ \ 
                   f(t,x,v) \ge C_1 \, e^{-C_2 \, |v|^K} 
 \end{equation} 
for some $C_1,C_2,K > 0$, does converge to equilibrium at rate ``almost exponential'', i.e. faster than 
any inverse power of $t$. 
This paper works in some {\it a priori} setting on 
the solution, since there is no general Cauchy theory whose solutions satisfies 
suitable estimates to apply the ``entropy-entropy production'' method. 
Nevertheless a natural question is wether the set of solution satisfying the {\it a priori} assumptions 
of~\cite{DV:eqEB:03} is not trivial, i.e. reduced to the spatially homogeneous case or to 
cases for which exponential convergence results are already known.  
Our study answers to this question, since 
the solutions in~\cite{Guo:mous:01} satisfy all the estimates of regularity and decay needed 
in~\cite{DV:eqEB:03}, and a consequence of Theorem~\ref{theo:main} below is that they 
also satisfy~\eqref{eq:hypEEP} (with $K=2$). 
\medskip

\subsection{Statement of the results} 
Now let us introduce the functional spaces and the macroscopic quantities on the solution. 
We define bounds uniform in space on the observables of the solutions. 
We shall study precisely in Section~\ref{sec:app} in which case there is a Cauchy 
theory which fits these assumptions. One can say briefly that they are satisfied 
at least for hard spheres and inverse power laws interactions, either  
in the spatially homogeneous setting, or in the spatially inhomogeneous setting 
for solutions ``in the small'' (i.e. for small time or near the equilibrium). 

Let us consider a function $f(t,x,v) \ge 0$ on $[0,T) \times \ens{T}^N \times \R^N$. We define 
its {\em local density} 
 \begin{equation*}
 \rho_f (t,x) := \int_{\R^N} f(t,x,v) \, dv,
 \end{equation*}
its {\em local energy}
 \begin{equation*}
 e_f (t,x) := \int_{\R^N} f(t,x,v) \, |v|^2 \, dv, 
 \end{equation*}
a weighted local energy 
 \begin{equation*}
 e_f ' (t,x) := \int_{\R^N} f(t,x,v) \, |v|^{\tilde{\gamma}} \, dv
 \end{equation*}
(where $\tilde{\gamma}$ is the 
positive part of $(2+\gamma)$),  
its {\em local entropy} 
 \begin{equation*}
 h_f (t,x) := - \int_{\R^N} f(t,x,v) \, \log f(t,x,v) \, dv,
 \end{equation*}
its local $L^p$ norm ($p \in [1,+\infty)$)
 \begin{equation*}
 l^p _f(t,x) := \| f (t,x,\cdot) \|_{L^p (\R^N _v)},   
 \end{equation*}
and its local $W^{2,\infty}$ norm~\footnote{see Subsection~\ref{subsec:not} for the notations}
 \begin{equation*}
 w _f (t,x) := \|f(t,x,\cdot)\|_{W^{2,\infty}  (\R^N _v)}.
 \end{equation*}
Note that in the sequel we shall systematically speak of hydrodynamical quantities on the 
solution in a generalized sense, since we include in this term the quantities $e_f '$, 
$h_f$, $l^p _f$ and $w _f$.

Then we define the following uniform bounds 
 \begin{equation*}
 \varrho_f := \inf_{(t,x) \in [0,T) \times \ens{T}^N} \rho_f (t,x), 
 \  \  \  \  \ 
 E_f := \sup_{(t,x) \in [0,T) \times \ens{T}^N} \left( e_f (t,x) + \rho_f (t,x) \right), 
 \end{equation*}
 \begin{equation*}
 E_f ' := \sup_{(t,x) \in [0,T) \times \ens{T}^N} e_f '(t,x), 
 \  \  \  \  \ 
 H_f := \sup_{(t,x) \in [0,T) \times \ens{T}^N} |h_f (t,x)|,
 \end{equation*}
 \begin{equation*}
 L_f ^p := \sup_{(t,x) \in [0,T) \times \ens{T}^N} l^p _f (t,x), 
 \ \ \ \ \ 
 W_f := \sup_{(t,x) \in [0,T) \times \ens{T}^N} w _f (t,x). 
 \end{equation*}

\Remark In the spatially homogeneous setting all these quantities are independent of the space 
variable $x$ and the uniformity in time is in most cases, well-known or obvious (see Section~\ref{sec:app}). 
\medskip

Our assumptions on the solution are as follows:
 \begin{itemize}
 \item When $\gamma \ge 0$ and $\nu < 0$ (hard or Maxwellian 
 potentials with cutoff) we shall assume that
  \begin{equation} \label{eq:hyp1} 
  \varrho_f > 0, \ \ \ E_f < + \infty, \ \ \ H_f < +\infty. 
  \end{equation}
 \item When $\gamma \in (-N,0)$ (singularity of the kinetic collision kernel) 
 we shall make the additional assumption that
  \begin{equation} \label{eq:hyp2} 
  L^{p_{\gamma}} _f < +\infty 
  \end{equation}
 with $p_{\gamma} > \frac{N}{N+\gamma}$ (notice that this uniform bound on $L^{p_{\gamma}} _f$ implies 
 the one on the local entropy). 
 \item When $\nu \in [0,2)$ (singularity of the angular collision kernel) we shall make 
 the additional assumption that
  \begin{equation} \label{eq:hypnc} 
  W_f < +\infty, \ \ \ E_f ' < +\infty
  \end{equation}
 (remark that when $\gamma \le 0$, we have $E_f ' \le E_f$ and the second part of this 
 assumption is not necessary). 
 \end{itemize}
\Remark Although the regularity part of the 
last assumption~\eqref{eq:hypnc} seems quite stronger compared to the other ones, 
the regularizing property of the non-cutoff collision operator often ensures that it holds, 
at least in some cases (see Section~\ref{sec:app}), and thus makes it rather natural. 
\medskip

We now state our main theorems. The first one deals with cutoff collision kernels. 
In this theorem a mild solution to the Boltzmann equation with initial datum $f_0$ is a 
function $f$ which satisfies~\eqref{eq:duham} pointwise (see Definition~\ref{def:sol} below).
 \begin{theorem} \label{theo:main}
 Let $B = \Phi \,b$ be a collision kernel which satisfies~\eqref{eq:B}, with 
 $\Phi$ satisfying~\eqref{eq:Phi} or~\eqref{eq:Phi:mol}, and 
 $b$ satisfying~\eqref{eq:b} with $\nu <0$. 
 Let $f(t,x,v)$ be a mild solution of the full Boltzmann equation in the torus  
 on some time interval $[0,T)$, $T \in (0,+\infty]$ such that 
  \begin{itemize}
  \item[(i)] if $\Phi$ satisfies~\eqref{eq:Phi} with $\gamma \ge 0$ or 
  if $\Phi$ satisfies~\eqref{eq:Phi:mol}, then $f$ satisfies~\eqref{eq:hyp1};
  \item[(ii)] if $\Phi$ satisfies~\eqref{eq:Phi} with $\gamma < 0$, then $f$ 
  satisfies~\eqref{eq:hyp1} and~\eqref{eq:hyp2}.
  \end{itemize}

 Then for all $\tau \in (0,T)$ there exists some $\rho >0$ and $\theta > 0$ depending 
 on $\tau$, $\varrho_f$, $E_f$, $H_f$ (and $L^{p_{\gamma}} _f$ if $\Phi$ 
 satisfies~\eqref{eq:Phi} with $\gamma < 0$), 
 such that for all $t \in [\tau,T)$ the solution is bounded from below by the uniform Maxwellian 
 distribution with density $\rho$ and temperature $\theta$, i.e.  
  \begin{equation}\label{eq:minor}
  \forall \, t \in [\tau,T), \ 
  \forall \, x \in \ens{T}^N, \ 
  \forall \, v \in \R^N, \ \ \ 
  f(t,x,v) \ge \rho \, \frac{e^{-\frac{|v|^2}{2 \theta}}}{(2\pi \theta)^{N/2}}.
  \end{equation}
 \end{theorem}
\Remarks Let us comment on the assumptions and conclusions of this theorem:
\smallskip 

1. The main case of application of this theorem one should think of is $B=|v-v_*|^\gamma \, b(\cos \theta)$ 
with $b$ bounded from above and below. It includes the hard spheres model (when $\gamma=1$ 
and $b=1$), and the so-called ``variable hard spheres'' model.
\smallskip

2. As the lower bound in~\eqref{eq:minor} does not depend 
on the space variable $x$, Theorem~\ref{theo:main} applies to spatially homogeneous solutions 
as well: take $f_0$ depending only on $v$ and $f(t,v)$ the corresponding solution of the 
homogeneous Boltzmann equation,  
then $f(t,v)$ is also 
a solution of the inhomogeneous Boltzmann equation in the torus and Theorem~\ref{theo:main} gives the 
appearance of a Maxwellian lower bound on the $v$ variable. 
This theorem thus includes and extends the previous 
result of Pulvirenti and Wennberg in~\cite{PW:lb:97}, 
valid for hard potentials. It gives new results for spatially homogeneous 
solutions in the case of soft potentials with cutoff (see Section~\ref{sec:app}).
\smallskip

3. In the inhomogeneous case, Theorem~\ref{theo:main} applies to the solutions of Ukai~\cite{Ukai:FBE:74} 
near the equilibrium for hard spheres, or to the solutions of Guo~\cite{Guo:mous:01} 
near the equilibrium for soft potentials with cutoff (see Section~\ref{sec:app}). 
\smallskip

4. More generally this theorem can be seen as an {\it a priori} result on the renormalized solutions 
(see~\cite{DiLi:89} and~\cite[Chapter~5]{CIP:94}), 
the only theory of solutions in the large at now. 
For instance for a gas of hard spheres in a torus, its converse says that if the solution $f$ vanishes, then 
either the local density $\rho_f$ has to vanish somewhere in the torus, or the local density $\rho_f$, 
energy $e_f$ or entropy $h_f$ have to blow-up somewhere in the torus. 
\medskip 

Now let us state the result we get for long-range interaction models, i.e.  
collision kernels with an angular singularity. 
In this theorem a mild solution to the Boltzmann equation with initial datum $f_0$ is a 
function $f$ which satisfies~\eqref{eq:duham:nc} pointwise (see Definition~\ref{def:sol:nc} below).
 \begin{theorem}\label{theo:main:nc}
 Let $B = \Phi \,b$ be a collision kernel which satisfies~\eqref{eq:B}, with 
 $\Phi$ satisfying~\eqref{eq:Phi} or~\eqref{eq:Phi:mol}, and 
 $b$ satisfying~\eqref{eq:b} with $\nu \in [0,2)$.   
 Let $f(t,x,v)$ be a mild solution of the full Boltzmann equation in the torus  
 on some time interval $[0,T)$, $T \in (0,+\infty]$ such that 
  \begin{itemize}
  \item[(i)] if $\Phi$ satisfies~\eqref{eq:Phi} with $\gamma \ge 0$ or 
  if $\Phi$ satisfies~\eqref{eq:Phi:mol}, then $f$ satisfies~\eqref{eq:hyp1} and~\eqref{eq:hypnc}; 
  \item[(ii)] if $\Phi$ satisfies~\eqref{eq:Phi} with $\gamma < 0$, then $f$ 
  satisfies~\eqref{eq:hyp1},~\eqref{eq:hyp2} and~\eqref{eq:hypnc}.
  \end{itemize}

 Then for all $\tau \in (0,T)$ and for any exponent $K$ such that 
  \[ K > 2 \, \frac{\log \left(2 + \frac{2 \nu}{2 - \nu}\right)}{\log 2}, \] 
 there exist $C_1 >0$ and $C_2 >0$ depending 
 on $\tau$, $K$, $\varrho_f$, $E_f$, $E' _f$, $H_f$, $W_f$ 
 (and $L^{p_{\gamma}} _f$ if $\Phi$ satisfies~\eqref{eq:Phi} 
 with $\gamma < 0$), such that
  \begin{equation*}
  \forall \, t \in [\tau,T), \ 
  \forall \, x \in \ens{T}^N, \ 
  \forall \, v \in \R^N, \ \ \ 
  f(t,x,v) \ge C_1 \, e^{-C_2 \, |v|^K}.
  \end{equation*}
 Moreover in the case when $\nu=0$, one can take $K=2$ (Maxwellian lower bound). 
 \end{theorem}

\Remarks Let us comment on the assumptions and conclusions of this theorem:
\smallskip 

1. One can check that this theorem is consistent with Theorem~\ref{theo:main} when 
$\nu \to 0$. Notice that the situation when $\nu=0$ is particular: the collision operator 
is non-cutoff and corresponds to some ``logarithmic derivative''.  
\smallskip 

2. This theorem is, to the knowledge of the author, the first quantitative lower bound 
result for non-cutoff collision kernels. It applies for instance to the spatially homogeneous solutions 
recently obtained in~\cite{DesvWenn:nc:pp} (see Section~\ref{sec:app}). 
\smallskip

3. We mention that an extension of Theorem~\ref{theo:main} and 
Theorem~\ref{theo:main:nc} to the case of a bounded open set $\Omega$ with 
specular reflection ou bounce-back boundary condition on $\partial \Omega$ is in progress, 
and will be treated in a second part of this work. 

\subsection{Methods of proof}

In the spatially homogeneous case the proof of~\cite{PW:lb:97} 
proceeds in two steps: first the construction of 
an ``upheaval point'' for the solution after a small time, i.e. a uniform minoration 
on a ball; second a ``spreading process'' of this bound from below after a small time by the 
collision process, iterated infinitely many times. 
Both these steps are based on a mixing property of the 
gain part of the collision operator, which is reminiscent of 
the regularization property of this gain part, discovered by Lions 
(see~\cite{Lion:rgainI+II:94,Lion:rgainIII:94} or~\cite[Chapter~2, Section~3.4]{Vill:hand} 
for a review).  
The second step was already present in~\cite{Carl:EB:32} and 
systematically used in~\cite{PW:lb:97}. The main contributions of our paper are: 
\begin{enumerate} 
\item a strategy to deal with space dependent solutions (Section~\ref{sec:proof}),  
based on an implementation of the ``upheaval'' and ``spreading'' steps along each caracteristic,  
keeping track carefully of the constants in order  
to get uniform estimates according to the choice of the characteristic;
\item a strategy to deal with non cutoff collision kernels (Section~\ref{sec:proof:nc}), 
based on the use of a suitable splitting of the collision operator between 
a cutoff part which still enjoys the spreading property, and a small non-cutoff part, 
for which we give $L^\infty$ estimates of smallness thanks to the regularity assumptions on the solution.  
A precise balance between these two parts then allows to obtain the lower bound in the 
non cutoff case, although slighlty weaker; 
\item the implementation of the general method for soft potentials as well (i.e for 
collision kernels with a singular kinetic part), and in any dimension (Sections~\ref{sec:proof} 
and~\ref{sec:proof:nc});
\item a detailed discussion of the connection between these results and the existing 
Cauchy theories (Section~\ref{sec:app}).  
\end{enumerate}
Here we adopt the point of view of an {\it a priori} setting which allows to treat separately the issue  
of the lower bound and the one of establishing {\it a priori} estimates on the solution. 
Therefore we do not adress the question of obtaining such {\it a priori} estimates in the general 
case, which is open at now. 
This point of view should be understood as a unified approach for all existing Cauchy theories, 
as well as a way to obtain {\it a priori} results when no Cauchy theory exists, or when the solutions are too weak.
\medskip

The paper runs as follows. Section~\ref{sec:tools} remains purely functional, Section~\ref{sec:proof} 
and~\ref{sec:proof:nc} work on arbitrary solutions in {\it a priori} setting, and only Section~\ref{sec:app} 
deals with solutions which have effectively been constructed by previous authors. 
Section~\ref{sec:proof} is devoted to the proof of Theorem~\ref{theo:main}, 
Section~\ref{sec:proof:nc} to the proof of Theorem~\ref{theo:main:nc} and Section~\ref{sec:app} 
applies these two theorems to the existing Cauchy theories. 

\subsection{Notation}\label{subsec:not}

In the sequel we shall denote $\langle \cdot \rangle = \sqrt{1 + |\cdot|^2}$. We define the 
weighted Lebesgue space $L^p _q (\R^N)$ ($p \in [1,+\infty]$, $q \in \R$) by the norm 
 \[ \| f \|_{L^p _q (\R^N)} = \left[ \int_{\R^N} |f (v)|^p \, \langle v \rangle^{pq} \, dv \right]^{1/p} \] 
if $p < +\infty$ and 
 \[ \| f \|_{L^\infty _q (\R^N)} = \sup_{v \in \R^N} |f (v)| \, \langle v \rangle^{q} \]
when $p = +\infty$. The Sobolev space $W^{k,p} (\R^N)$ ($p \in [1,+\infty]$ and $k \in \N$) 
is defined by 
 \[ \| f \|_{W^{k,p} (\R^N)} = 
       \left[ \sum_{|s| \le k} \|\partial^s f (v)\|^p _{L^p} \right]^{1/p}, \]
with the usual notation $H^k = W^{k,2}$. 
Concerning the collision kernel we define the $L^1$ norm of $b$ on the unit sphere 
when $\nu < 0$ (integrable angular collision kernel) by 
 \[ n_b = \int_{\ens{S}^{N-1}} b(\cos \theta) \, d\sigma = |\ens{S}^{N-2}| \, 
          \int_0 ^\pi b(\cos \theta) \, \sin^{N-2} \theta \, d\theta, \]
and in the case $\nu \in [0,2)$ we define 
 \[ m_b = \int_{\ens{S}^{N-1}} b(\cos \theta) \, (1 - \cos \theta) \, d\sigma = |\ens{S}^{N-2}| \, 
          \int_0 ^\pi b(\cos \theta) \, (1 - \cos \theta) \, \sin^{N-2} \theta \, d\theta, \]
which is always finite (since $\nu < 2$), 
and is related to the cross-section for momentum transfer (see~\cite[Chapter~1, Section~3.4]{Vill:hand}). 
Finally we define
 \[ \ell_b = \inf_{\pi/4 \le \theta \le 3 \pi/4} b(\theta) \] 
which is strictly positive by assumption. 

In the following we shall keep track explicitely of the dependency of the constants according to 
the bounds on the collision kernel and the estimates on the solution. As a convention, 
``$\mbox{cst}$'' shall systematically denote any constant depending only on the dimension $N$, 
$\gamma$, $\nu$ and $b_0$. For a real $x$, we shall denote $x^+$ the positive part of $x$ and we 
recall the shorthand $\tilde{\gamma} = (2+\gamma)^+$. 

\section{Functional toolbox}\label{sec:tools}
\setcounter{equation}{0}

In this section we shall gather functional tools useful for the sequel. 
On one hand, Lemma~\ref{lem:L}, Lemma~\ref{lem:Q1} are precised  
versions of well-known results adapted to our study: we need $L^\infty$ estimates 
whereas the usual framework of such estimates are integral spaces. 
On the other hand, Lemma~\ref{lem:up} and Lemma~\ref{lem:spread} are essentially 
generalizations of results in~\cite{PW:lb:97}. 
We extend these results for any cutoff potentials 
(in the sense of~\eqref{eq:Phi} and~\eqref{eq:b} with $\nu < 0$), 
in any dimension. Moreover we intend to use these results in the 
context of spatially inhomogeneous solutions, using the fact that the collision operator 
is local in $t$ and $x$, which allows to see these variables as parameters in the functional 
estimates. Thus we shall track precisely the dependence of these estimates on the hydrodynamical 
quantities: $\rho_f$, $e_f$, $e' _f$, $h_f$, $l^p _f$, $w_f$.

\subsection{The cutoff case}

We introduce Grad's splitting 
 \begin{eqnarray*}
 Q (g,f) &=& Q^+ (g,f) - Q^- (g,f) \\
 Q^+(g,f) &:=& \int _{\R^N} \, dv_*
 \int _{\ens{S}^{N-1}} \, d\sigma B(|v-v_*|, \cos \theta) \, g'_* f' \\
 Q^-(g,f) &:=& \int _{\R^N} \, dv_*
 \int _{\ens{S}^{N-1}} \, d\sigma B(|v-v_*|, \cos \theta) \, g_* f
 \end{eqnarray*}
where $Q^+$ is called the gain term and $Q^-$ is called the loss term. We write the 
loss term as 
 \[ Q^-(g,f) = L[g] \, f \] 
with 
 \begin{equation}\label{eq:defL}
 L[g(t,x,\cdot)](v) = n_b \, \int_{\R^N} \Phi(v-v_*) \, g(t,x,v_*) \, dv_*. 
 \end{equation}

First let us give an $L^\infty$ bound on the action of the loss term along the characteristics. 
 \begin{lemma}\label{lem:L}
 Let $g$ be a mesurable function on $\R^N$. Then 
  \begin{itemize}
  \item[(i)] If $\Phi$ satisfy~\eqref{eq:Phi} with $\gamma \ge 0$ or 
  if $\Phi$ satisfies~\eqref{eq:Phi:mol}, then 
   \begin{equation*}
   \forall \, v \in \R^N, \ \ \ |g * \Phi (v)| \le 
   \mbox{{\em cst}} \; C_\Phi \, \|g\|_{L^1 _2} \, \langle v \rangle^{\gamma ^+}.  
   \end{equation*}
  \item[(ii)] If $\Phi$ satisfy~\eqref{eq:Phi} with $\gamma < 0$, then 
   \begin{equation*}
   \forall \, v \in \R^N, \ \ \ |g * \Phi (v)| \le \mbox{{\em cst}} \; C_\Phi \,  
   \left[ \|g\|_{L^1 _2} + \|g\|_{L^p} \right] \, \langle v \rangle^{\gamma ^+}. 
   \end{equation*}
  with $p > N/(N+\gamma)$.
  \end{itemize}
 \end{lemma}

 \begin{corollary}\label{coro:LS}
 As a straightforward consequence we obtain the following 
 estimates on the operators $L$ and $S$ defined respectively in~\eqref{eq:defL} 
 and~\eqref{eq:defS}
  \begin{equation}\label{eq:LS} 
  \forall \, v \in \R^N, \ \ \ |L[g](v)| \le C_L \, \langle v \rangle^{\gamma ^+} \ \  
  \mbox{ and } \ \ |S[g](v)| \le C_S \, \langle v \rangle^{\gamma ^+}
  \end{equation}
 where the constants $C_L$ and $C_S$ are defined by:
  \begin{itemize} 
  \item[(i)] If $\Phi$ satisfy~\eqref{eq:Phi} with $\gamma \ge 0$ or 
  if $\Phi$ satisfies~\eqref{eq:Phi:mol}, then 
   \begin{equation*} 
   C_L = \mbox{{\em cst}} \; n_b \, C_\Phi \, e_g, \ \ \ \ \ C_S = \mbox{{\em cst}} \; m_b \, C_\Phi \, e_g.
   \end{equation*}
  \item[(ii)] If $\Phi$ satisfy~\eqref{eq:Phi} with $\gamma < 0$, then 
   \begin{equation*} 
   C_L = \mbox{{\em cst}} \; n_b \, C_\Phi \, \left[ e_g + l^p _g \right], \ \ \ \ \ 
   C_S = \mbox{{\em cst}} \; m_b \, C_\Phi \, \left[ e_g + l^p _g \right].
   \end{equation*}
  \end{itemize}
 \end{corollary}

\begin{proof}[Proof of Lemma~\ref{lem:L}]
In the case $\Phi$ satisfies~\eqref{eq:Phi} with $\gamma \ge 0$ or~\eqref{eq:Phi:mol}, 
the proof is obvious and amounts to a triangular inequality. 

In the case $\Phi$ satisfies~\eqref{eq:Phi} with $\gamma \in (-N,0)$, one can split $g * \Phi (v)$ into 
 \begin{multline*}
 g * \Phi (v) = \int_{\{ v_* \ ; \ |v-v_*| \le 1 \}} \Phi(v-v_*) \, g(v_*) \, dv_* \\ 
                + \int_{\{ v_* \ ; \ |v-v_*| \ge 1 \}} \Phi(v-v_*) \, g(v_*) \, dv_*
              =: I_1 + I_2.
 \end{multline*}
On one hand, 
 \[ \forall \, v \in \R^N, \ \ \ |I_2 (v)| \le C_\Phi \, \|g\|_{L^1} \le C_\Phi \, \|g\|_{L^1 _2} \]
and on the other hand, by Cauchy-Schwarz inequality
 \[ \forall \, v \in \R^N, \ \ \ |I_1 (v)| \le C_\Phi \, 
    \left[ \int_{\{ v_* \ ; \ |v-v_*| \le 1 \}} |v-v_*|^{\gamma p'} \, dv_* \right]^{1/p'} \, 
    \| g \|_{L^p} \] 
which gives the result since 
 \[ \mbox{cst} = \int_{\{ v_* \ ; \ |v-v_*| \le 1 \}} |v-v_*|^{\gamma p'} \, dv_* = 
    \int_{\{ u \ ; \ |u| \le 1 \}} |u|^{\gamma p'} \, du < +\infty \] 
as soon as $\gamma p' > -N$, i.e $p>\frac{N}{N+\gamma}$. 
\end{proof}

The next lemma uses the mixing property of $Q^+$ in order to obtain a minoration 
of $Q^+ ( Q^+ (\cdot, \cdot), \cdot)$ on a ball. This will be the starting 
point for the construction of an ``upheaval point'' by the iterated Duhamel 
formula. 
 \begin{lemma}\label{lem:up}
 Let $B = \Phi \,b$ be a collision kernel which satisfies~\eqref{eq:B}, with 
 $\Phi$ satisfying~\eqref{eq:Phi} or~\eqref{eq:Phi:mol}, and 
 $b$ satisfying~\eqref{eq:b} with $\nu \le 0$. Let $g(v)$ be a nonnegative function on $\R^N$  
 with bounded energy $e_g$ and entropy $h_g$ and a mass $\rho_g$ such that $0<\rho_g<+\infty$.  
 Then there exist $R_0 > 0$, $\delta_0 >0$, $\eta_0 >0$ and $\bar{v} \in B(0,R_0)$ such that 
  \begin{equation*}\label{eq:up}
  Q^+ ( Q^+ ( g \, {\bf 1}_{B(0,R_0)} , g \, {\bf 1}_{B(0,R_0)}), g \, {\bf 1}_{B(0,R_0)}) 
  \ge \eta_0 \, {\bf 1}_{B(\bar{v},\delta_0)} 
  \end{equation*}
 and $R_0$, $\delta_0$, $\eta_0$ only depend on $B$, on a lower bound on $\rho_g$, 
 and upper bounds on $e_g$ and $h_g$. 
 \end{lemma}
\Remark Another strategy to obtain this ``upheaval point'' for hard potentials 
could have been to iterate the regularity property of $Q^+$ in the form proved 
in~\cite{MV:RegHom:03} in Sobolev spaces enough times to get some continuous function.  
We did not follow this method which is less direct, and 
leads to harder computations to track the exact dependence of the constant. Nevertheless 
the remark emphasizes the fact that the mixing property used here on $Q^+$ is a {\em linear} one,  
which is consistent with the regularity theory for this operator 
(see the regularity theory of $Q^+$ in the cutoff case in~\cite{MV:RegHom:03}).
\medskip

\begin{proof}[Proof of Lemma~\ref{lem:up}]
This lemma is a slight variant of~\cite[Lemma~3.1]{PW:lb:97}, whose proof can be straightforwardly 
adapted here. Note that this proof was made assuming that $b$ is bounded below by a positive 
quantity on the whole interval $[0,\pi]$ but as pointed out in~\cite[Proof of Lemma~3.1]{PW:lb:97} 
the proof still works the same under the sole assumption on $b$ that it is bounded 
below by a positive quantity near $\theta \sim \pi/2$, which is satisfied under our assumptions on $b$. 
 
Therefore we only generalize the formula in the proof to any dimension and to any power $\gamma$ of 
the collision kernel, and we precise the dependence of $R_0$, $\delta_0$ and $\eta_0$ 
according to the quantities $\rho_g$, $e_g$ and $h_g$.

First let us suppose that $\Phi$ satisfies~\eqref{eq:Phi} with $\gamma \ge 0$, in order 
to satisfy the assumptions of~\cite[Lemma~3.1]{PW:lb:97}.  
As for $R_0$, in the proof of~\cite[Lemma~3.1]{PW:lb:97} $R_0$ is chosen such that 
 \[ \int_{|v|\le R_0} g(v) \, dv \ge \frac{\rho_g}2. \] 
The estimate
 \[ \int_{|v|\ge R_0} g(v) \, dv \le \frac{e_g}{R_0 ^2} \] 
yields the possible choice 
 \[ R_0 = \sqrt{\frac{2 \, e_g}{\rho_g}}. \]
Then it is straigthforward to see that $\delta_0$ depends only on upper bounds on $e_g$ and 
$h_g$, and 
 \[ \eta_0 = \mbox{cst} \, \ell_b \, c_\Phi \, R_0 ^{\gamma - (3N-1)} \, \delta_0^{2N}. \]

The case of a mollified kinetic collision kernel $\Phi$ (assumption~\eqref{eq:Phi:mol}) with 
$\gamma \ge 0$ reduces to the case of a kinetic collision kernel satisfying~\eqref{eq:Phi} 
with $\gamma \ge 0$. 
Indeed when $\gamma \ge 0$, we have the bound from below $\Phi(z) \ge c_\Phi \, |z|^\gamma$ for all 
$z \in \R^N$ and the proof is unchanged. 

When $\Phi$ satisfies~\eqref{eq:Phi} or~\eqref{eq:Phi:mol} with $\gamma < 0$, 
we first choose $R_0$ as above, with $R_0 \ge 1$ (it is possible up to take a 
bigger $R_0$). Then we use that on $B(0,R_0)$, we have that 
$\Phi(z) \ge c_\Phi \, R_0 ^{\gamma}$, which means we can apply the 
formula above for the case $\gamma = 0$, and the final formula for $\eta_0$ is 
unchanged. 
\end{proof}

\begin{figure}[h]
 \epsfysize=6cm
 $$\epsfbox{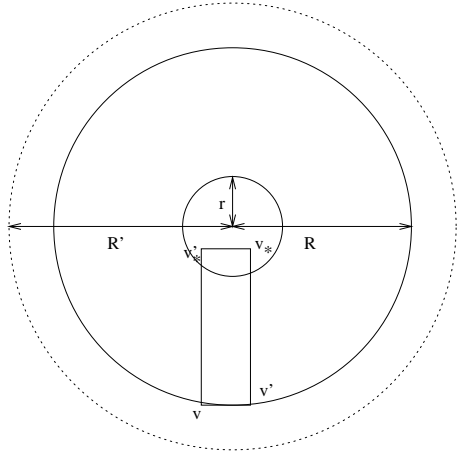}$$
 \caption{Spreading property of $Q^+$}\label{fig:col}
 \end{figure}

The next lemma gives a precise estimate of the ``spreading property'' of $Q^+$ 
(according to the velocity variable), which is pictured in Figure~\ref{fig:col}: 
for any $R' < \sqrt{r^2+R^2}$, for any $v$ in the ball with radius $R'$, it is possible 
to find collisions with post-collision velocity $v,v_*$ and 
taking the pre-collision velocity $v'_*$ inside the ball with radius $R$ and 
the pre-collision velocity $v'$ inside the ball with radius $r$.   

 \begin{lemma}\label{lem:spread}
 Let $B = \Phi \,b$ be a collision kernel which satisfies~\eqref{eq:B}, with 
 $\Phi$ satisfying~\eqref{eq:Phi} or~\eqref{eq:Phi:mol}, and 
 $b$ satisfying~\eqref{eq:b} with $\nu \le 0$. 
 Then for any $\bar{v} \in \R^N$, $0 < r \le R$, $\xi \in (0,1)$, we have 
  \begin{equation}\label{eq:spread}
  Q^+ ( {\bf 1}_{B(\bar{v},R)} , {\bf 1}_{B(\bar{v},r)} )
  \ge \mbox{{\em cst}} \; \ell_b \, c_\Phi \, r^{N-3} \, R^{3+\gamma} \, \xi^{N/2-1} \, 
  {\bf 1}_{B\big(\bar{v},\, \sqrt{r^2 + R^2} (1-\xi)\big)}.
  \end{equation}

 As a consequence in the particular quadratic case $\delta=r=R$, we obtain
  \begin{equation}\label{eq:spread:quad}
  Q^+ ( {\bf 1}_{B(\bar{v},\, \delta)} , {\bf 1}_{B(\bar{v},\, \delta)}) 
  \ge \mbox{{\em cst}} \; \ell_b \, c_\Phi \, \delta^{N+\gamma} \, \xi^{N/2-1} \, 
  {\bf 1}_{B\big(\bar{v},\, \delta \sqrt{2} (1-\xi)\big)}
  \end{equation}
 for any $\bar{v} \in \R^N$ and $\xi \in (0,1)$.
 \end{lemma}

\Remark In the sequel we shall use the quadratic version of Lemma~\ref{lem:spread} 
(i.e when $r=R$) which seems compulsory when one wants to obtain the 
optimal Maxwellian decrease at infinity for the lower bound. 
Nervertheless we give a bilinear version since it highlights again the fact 
that the ``spreading effect'' of $Q^+$ is a linear one.
\medskip

\begin{proof}[Proof of Lemma~\ref{lem:spread}]
This result is a bilinear version of~\cite[Lemma~3.2]{PW:lb:97}, written here in any 
dimension and for any power $\gamma$. Thus we only recall the main steps of the proof, 
especially those where the bilinearity, the dimension $N$, or $\gamma$ play some role.

First we deal with collision kernel such that $\Phi$ satisfies~\eqref{eq:Phi} with 
$\gamma \ge 0$. As a general property $Q^+$ satisfies the homogeneity relation
 \[ Q^+ (g,f) (\lambda v) = \lambda^{N+\gamma} \, Q^+ ( g(\lambda \cdot), f(\lambda \cdot))(v) \] 
and the invariance by translation allows to reduce the proof of~\eqref{eq:spread} to the proof of 
  \begin{equation*}\label{eq:spread:bis}
  Q^+ ( {\bf 1}_{B(0,\, 1)} , {\bf 1}_{B(0,\, p)})
  \ge \mbox{cst} \, \ell_b \, c_\Phi \, p^{N-3} \, \xi^{N/2-1} \, 
  {\bf 1}_{B\big(0,\sqrt{1 + p^2} (1-\xi)\big)}
  \end{equation*}
with $p \le 1$ stands for $r/R$.  
Now by isotropic invariance we can assume $v=z \, {\bf e_N}$ 
with $z < \sqrt{1 + p^2}$ and $({\bf e_1}, \dots, {\bf e_N})$ an orthonormal basis. 
By Carleman representation~\cite{Carl:57}, we have 
  \begin{multline*}
  Q^+ ( {\bf 1}_{B(0,\, 1)} , {\bf 1}_{B(0,\, p)}) (v= z \, {\bf e_N})
  \ge \\ c_\Phi \, \int_{v' \in \R^N} \frac{{\bf 1}_{B(0,\, p)} (v')}{|v-v'|^{N-1-\gamma}} 
      \left[ \int_{v' _* \in E_{vv'}} {\bf 1}_{B(0,\, 1)} (v' _*) 
      \tilde{b}(\theta) \, dv' _* \right] \, dv' 
  \end{multline*}
where $E_{vv'}$ is the hyperplan orthogonal to $v'-v$ and containing $v$. 
We write the integral along $v'$ in spherical coordinates centered in $v$ 
and we use the bound from below 
$\tilde{b}(\theta) = b(\theta) \, (\sin \theta/2)^{-\gamma} \ge \mbox{cst} \; \ell_b$ 
for $\theta \in [\pi/4,3\pi/4]$ given by the assumptions on $b$
  \begin{multline*}
  Q^+ ( {\bf 1}_{B(0,\, 1)} , {\bf 1}_{B(0,\, p)}) (v= z \, {\bf e_N})
  \ge \\ \mbox{cst} \; \ell_b \, c_\Phi \, 
      \int_0 ^{+\infty} \int_{\ens{S}^{N-1}} {\bf 1}_{B(0,\, r)} (v+\rho \sigma) \, 
      \rho^\gamma \, \mbox{Vol }( E_{vv'} \cap B(0,\, 1) \cap \cal{C}_{v,\rho}) \, d\rho \, d\sigma,  
  \end{multline*}
where 
 \[ \cal{C}_{v,\rho} = \left\{ u \in \R^N, \ \ 
    \tan \frac{\pi}{8} \, \rho \le |u-v| \le \tan \frac{3\pi}{8} \, \rho \right\}. \]
Finally it is easy to see that the integrand is invariant under rotation around 
the axis $(0, {\bf e_N})$, which allows to simplify the part of integration over the 
unit sphere $\ens{S}^{N-2}$ of the hyperplan orthogonal to $(0, {\bf e_N})$:
  \begin{multline*}
  Q^+ ( {\bf 1}_{B(0,\, 1)} , {\bf 1}_{B(0,\, p)}) (v= z \, {\bf e_N})
  \ge \\ \mbox{cst} \; \ell_b \, c_\Phi \, \int_0 ^{+\infty} \int_0 ^\pi  
  (\sin \alpha)^{N-2} \, {\bf 1}_{B(0,\, r)} (v+\rho \sigma) \, 
  \rho^\gamma \, \mbox{Vol }( E_{vv'} \cap B(0,\, 1) \cap \cal{C}_{v,\rho}) \, d\rho \, d\alpha. 
  \end{multline*}
where $\sigma$ in the integrand stands for any $-(\cos \alpha){\bf e_N} + (\sin \alpha) {\bf u}$  
with ${\bf u}$ is any vector of the set of unit vectors orthogonal to ${\bf e_N}$. 
Some elementary geometrical computation lead to 
 \[ \mbox{Vol }( E_{vv'} \cap B(0,\, 1) ) = 
    \mbox{cst} \, \left( 1 - z^2 \cos^2 \alpha \right)^{\frac{N-1}2} \, 
    {\bf 1}_{\{\cos^2 \alpha \le 1/z^2 \}} \]
and shows that 
 \[ E_{vv'} \cap B(0,\, 1) \subset \cal{C}_{v,\rho} \] 
when 
 \[ a \, \left( z \sin \alpha + \sqrt{1 - z^2 \cos^2 \alpha} \right) 
 \le \rho \le b \, \left( z \sin \alpha - \sqrt{1 - z^2 \cos^2 \alpha} \right)\]
with 
 \[ a := \tan^{-1} \frac{3\pi}{8} < 1, \ \ \ b:= \tan^{-1} \frac{\pi}{8} > 1. \]
This inequality is possible as soon as 
 \[ 1 - z^2 \cos^2 \alpha \le \lambda^2 \, z^2 \sin^2 \alpha \]
with $\lambda = (b - a)/(b + a) < 1$. 
If one sets $y=z \cos \alpha$ as a new variable, this inequality means 
 \[ y \ge \sqrt{\frac{1 - \lambda^2 \, z}{1-\lambda^2}} \] 
and one gets 
 \begin{multline*}
 Q^+ ( {\bf 1}_{B(0,\, 1)} , {\bf 1}_{B(0,\, p)}) (v= z \, {\bf e_N})
 \ge \\ \mbox{cst} \; \ell_b \, c_\Phi \, 
 \int_{\max\left\{\sqrt{z^2-1},\sqrt{\frac{1 - \lambda^2 \, z^2}{1-\lambda^2}}\right\}} ^1 
 \int_{\max\{y-\sqrt{1-z^2+y^2},\, a (\sqrt{z^2-y^2}+\sqrt{1-y^2})\}} 
      ^{\min\{y+\sqrt{1-z^2+y^2},\, b (\sqrt{z^2-y^2}-\sqrt{1-y^2})\}}  \\  
 \rho^\gamma \, \left( 1 - \frac{y^2}{z^2} \right)^{\frac{N-3}2} 
 \, \left( 1 - y^2 \right)^{\frac{N-1}2} \, dy \, d\rho. 
 \end{multline*}
Now setting $z = \sqrt{1 + p^2}(1-\xi)$ and computing an expansion of this expression 
according to $\xi$ in the same way as in the end of the proof of~\cite[Lemma~3.2]{PW:lb:97}, 
one gets the following estimates:
 \[ \left( 1 - y^2 \right)^{\frac{N-1}2} = \big( 2(1+p^2) \xi \big)^{\frac{N-1}2} 
                                           + O(\xi^{\frac{N+1}2}) \] 
 \[ \left( z^2 - y^2 \right)^{\frac{N-3}2} = p^{N-3} + O(\xi) \] 
 \[ \int_{\max\{y-\sqrt{1-z^2+y^2},\, a (\sqrt{z^2-y^2}+\sqrt{1-y^2})\}} 
        ^{\min\{y+\sqrt{1-z^2+y^2},\, b (\sqrt{z^2-y^2}-\sqrt{1-y^2})\}} \rho^\gamma \, d\rho 
    = \sqrt{2(1+p^2)\xi -2(1-y)} + O(\xi^{3/2}). \]
Then similar computations as in the proof of~\cite[Lemma~3.2]{PW:lb:97} conclude the proof 
(for the integration on $y$, the condition $y \ge \sqrt{(1 - \lambda^2 \, z^2)/(1-\lambda^2)}$ 
plays no role at the limit since 
$\sqrt{(1 - \lambda^2 \, z^2)/(1-\lambda^2)} \to_{\xi \to 0} \mbox{cst} <1$). 

As for the previous lemma the case of a mollified kinetic collision kernel 
$\Phi$ (assumption~\eqref{eq:Phi:mol}) with 
$\gamma \ge 0$ reduces the case of a kinetic collision kernel satisfying~\eqref{eq:Phi} 
with $\gamma \ge 0$ since we have the bound from below $\Phi(z) \ge c_\Phi \, |z|^\gamma$ for all 
$z \in \R^N$ and the proof is unchanged. 

When $\Phi$ satisfies~\eqref{eq:Phi} or~\eqref{eq:Phi:mol} with $\gamma < 0$, 
we use that on $B(0,R)$, we have that 
$\Phi(z) \ge c_\Phi \, R ^{\gamma}$ (assuming $R \ge 1$ without restriction for the 
sequel) which means we can apply the 
formula above for the case $\gamma = 0$, and the final formula~\eqref{eq:spread} is 
unchanged. 
\end{proof}

\subsection{The non-cutoff case}

The two lemmas below will be useful in the treatment of non-cutoff collisions kernels. 
They express the fact that non-grazing collisions constitute the dominant 
term of the collision operator as long as ``spreading effect'' is concerned. They are essentially 
based on the by now well-known idea of using symmetry-induced cancellations effects 
in order to deal with the angular singularity (see~\cite{Vill:new:98,ADVW:00,AlexVill:nc:02}). 

In the case of non cutoff collision kernels, the usual Grad's splitting 
$Q=Q^+ - Q^-$ does not make sense anymore. However, the following splitting still makes sense:
 \begin{eqnarray} \nonumber \label{eq:split:nc}
 Q(g,f) &=& \int_{\R^N} dv_* \int_{\ens{S}^{N-1}} d\sigma \, B \, g'_* \, (f'-f) + 
            \left( \int_{\R^N} dv_* \int_{\ens{S}^{N-1}} d\sigma \, B \, (g'_* - g_*) \right) f \\
        &=:& Q^1 + Q^2 
 \end{eqnarray}
Thanks to the cancellation Lemma~\cite[Lemma~1]{ADVW:00}, the operator $Q^2$ can be written
 \[ Q^2(g,f) = S[g] \, f \] 
with 
 \begin{equation}\label{eq:defS}
 S[g] (v) := |\ens{S}^{N-2}| \, \left( \int_0 ^{\pi/2} \sin^{N-2} \theta \, 
 \left[ \frac{1}{\cos^{N+\gamma}(\theta/2)} - 1 \right] \, b(\theta) \, d\theta \right)\, \Phi * g (v). 
 \end{equation}

Corollary~\ref{coro:LS} gives $L^\infty$ estimates~\eqref{eq:LS} on $S$. 
Now let us turn to the $L^\infty$ estimates for $Q^1$. 
 \begin{lemma}\label{lem:Q1}
 Let $B = \Phi \,b$ be a collision kernel which satisfies~\eqref{eq:B}, with 
 $\Phi$ satisfying~\eqref{eq:Phi} or~\eqref{eq:Phi:mol}, and 
 $b$ satisfying~\eqref{eq:b} with $\nu \in [0,2)$. 
 Let $f,g$ be mesurable functions on $\R^N$. Then 
  \begin{itemize}
  \item[(i)] If $\Phi$ satisfy~\eqref{eq:Phi} with $2+\gamma \ge 0$ or 
  if $\Phi$ satisfies~\eqref{eq:Phi:mol}, then 
   \begin{equation*}\label{eq:Q1}
   \forall \, v \in \R^N, \ \ \ |Q^1 (g,f) (v)| \le 
   \mbox{{\em cst}} \; m_b \, C_\Phi \, 
   \|g\|_{L^1 _{\tilde{\gamma}}} \, \|f\|_{W^{2,\infty}} \, \langle v \rangle^{\tilde{\gamma}}. 
   \end{equation*}
  \item[(ii)] If $\Phi$ satisfy~\eqref{eq:Phi} with $2+\gamma < 0$, then 
   \begin{equation*}\label{eq:Q1:lp}
   \forall \, v \in \R^N, \ \ \ |Q^1 (g,f) (v)| \le 
   \mbox{{\em cst}} \; m_b \, C_\Phi \, \left[ \|g\|_{L^1 _{\tilde{\gamma}}} + \|g\|_{L^p} \right]\, 
   \|f\|_{W^{2,\infty}} \, \langle v \rangle^{\tilde{\gamma}}.
   \end{equation*}
  with $p > N/(N+\gamma+2)$.
  \end{itemize}
 \end{lemma} 
\Remarks 
1. In the treatment of $Q^2$, the ``derivative-like'' difference
$(f'-f)$ can be transferred to the angular part of the collision kernel by a process of change of variable which 
plays the same role as integration by part for classical differential operators. 
One can not do the same in Lemma~\ref{lem:Q1} because there is no ``decoupling'' 
of the two arguments $f$ and $g$. 
\smallskip

2. The proof of Lemma~\ref{lem:Q1} is based on a similar idea as cancellation 
lemmas in~\cite{Vill:new:98,ADVW:00,AlexVill:nc:02}. We make a Taylor expansion of 
$(f'-f)$ for small deviation angles in order to compensate for 
the grazing collision singularity of the collision kernel. However for smooth functions 
the quantity $(f'-f)$ compensate only for a singularity of order $1$. 
Thus one has to use the symmetry of the collision sphere in 
order to compensate for strong singularities: 
to compensate for a singularity of order $2$, one has to kill the first order term 
of the Taylor expansion of $(f'-f)$ by integrating on the $(N-2)$-dimensional sphere included 
in the collision sphere which is orthogonal to the vector $v-v_*$ and contains $v'$. 
\smallskip

3. Note that Lemma~\ref{lem:Q1} is reminiscent of~\cite[Proposition~4]{AlexVill:nc:02} 
which implies that the complete non cutoff operator $Q$ satisfies the following inequality 
 \[ \|Q(g,f) \|_{W^{-2,1}} \le C \, \|g\|_{L^1 _{\tilde{\gamma}}} \, \|f\|_{L^1 _{\tilde{\gamma}}} \] 
Essentially the difference is that Lemma~\ref{lem:Q1} is 
intended to provide an $L^\infty$ control. This is why it requires an $L^p$ bound on the solution 
for soft potentials.  
The proof of~\cite[Proposition~4]{AlexVill:nc:02} uses a duality argument 
and the pre-post collisional change of variable to pass the ``derivative-like'' difference on the dual 
test function. The $(N-2)$-dimensional sphere on which cancellations occur does not appear explicitely 
in the representation formula as in our proof, but is rather implicit in a projection argument. 
Here we proceed directly, using Carleman representation. 
\medskip

\begin{proof}[Proof of Lemma~\ref{lem:Q1}]
In order to isolate the exact sphere on which we want to use symmetry properties, 
we use the Carleman representation (exchanging the roles of $v'$ and $v'_*$) 
 \begin{equation*}
 Q^1 (g,f)(v) = \int_{\R^N} dv' _* \, \frac{g' _*}{|v-v'_*|^{N-1}}  \int_{E_{v,v'_*}} dv' \, 
                b(\cos \theta) \, \Phi(v-v_*) \,(f'-f)
 \end{equation*}
where $E_{v,v'_*}$ denotes the hyperplan orthogonal to the vector $v-v'_*$ and 
containing $v$. Now let us write the integration of the $v'$ variable in spherical 
coordinate of center $v$, i.e $v' = v + \rho \sigma$ where $\rho \in \R_+$ and 
$\sigma$ describes the $(N-2)$-dimensional unit sphere of $E_{v,v'_*}$, 
denoted by $\ens{S}_{v-v'_*} ^{N-2}$.  
 \begin{equation*}
 Q^1 (g,f)(v) = \int_{\R^N} dv' _* \, \frac{g' _*}{|v-v'_*|^{N-1}}  \int_0 ^\infty 
                d\rho \, b(\cos \theta) \, \Phi(v-v_*) \, \rho^{N-2} 
                \, \Bigg( \int_{\ens{S}_{v-v'_*} ^{N-2}} d\sigma \, 
                \big(f(v+ \rho \sigma)-f(v) \big) \Bigg).
 \end{equation*}
Now let us study more precisely the quantity 
 \begin{equation*}
 I = \int_{\ens{S}_{v-v'_*} ^{N-2}} d\sigma \, \big(f(v+ \rho \sigma)-f(v) \big).
 \end{equation*}
If $\nabla f$ denotes the gradient of $f$ and $\nabla^2 f$ its Hessian 
matrix, one has the following Taylor expansion: 
 \begin{equation*}
 f(v+ \rho \sigma) = f(v) + \rho \big( \nabla f(v) \cdot \sigma \big)
                     + \frac{\rho^2}{2} \, \left\langle \nabla^2 f (v + \rho' \sigma) 
                     \cdot \sigma , \sigma \right\rangle
 \end{equation*}
where $0 \le \rho' \le \rho$. By bounding the last term and 
taking the integral over $\sigma$, we get the estimate
 \begin{equation*}
 \left| I - \rho \, \left( \int_{\ens{S}_{v-v'_*} ^{N-2}} d\sigma \, 
     \big(\nabla f(v) \cdot \sigma \big) \right) \right| 
 \le  \frac{\rho^2}{2} \, |\ens{S}^{N-2}| \, \|f\|_{W^{2,\infty}}.
 \end{equation*}
As the term involving $\nabla f$ vanishes by symmetry, we obtain 
 \begin{equation*}
 |I| \le \frac{\rho^2}{2} \, |\ens{S}^{N-2}| \, \|f\|_{W^{2,\infty}}.
 \end{equation*}
Thus we get for some function $\phi$ in $L^1 _{\tilde{\gamma}}$ 
 \begin{multline*}
 \left|\int_{\R^N} Q^1 (g,f)(v) \phi(v) \, dv \right| \\ \le 
    \|f\|_{W^{2,\infty}} \, |\ens{S}^{N-2}| \, \int_{\R^N} dv \int_{\R^N} dv' _* \, 
                \frac{g' _*}{|v-v'_*|^{N-1}}  \int_0 ^\infty 
                d\rho \, b(\cos \theta) \, \Phi(v-v_*)  
                \, \frac{\rho^N}{2} \, |\phi(v)| \\ \le
    \|f\|_{W^{2,\infty}} \,  \int_{\R^N} dv \int_{\R^N} dv' _* \, \frac{g' _*}{|v-v'_*|^{N-1}}  \int_0 ^\infty 
                d\rho \, \rho^{N-2} \int_{\ens{S}_{v-v'_*} ^{N-2}} d\sigma \, 
               b(\cos \theta) \, \Phi(v-v_*) \, \frac{\rho^2}{2} \, |\phi(v)| \\ \le
    C_\Phi \, \|f\|_{W^{2,\infty}} \, 
                    \int_{\R^{2N} \times \ens{S}^{N-1}} 
                    dv \, dv_* \, d\sigma \, b(\cos \theta) \, 
                    \big( \sin (\theta/2) \big)^2 \, |v-v_*|^{2+\gamma} \, 
                    |g'_*| \, |\phi|. 
 \end{multline*}
Finally we cut the integral in two parts, for $\theta \in [0,\pi/2]$ and 
for $\theta \in [\pi/2,\pi]$. 
For the first part we use the pre-postcollisional change of variable and the change of 
variable $(v,v_*,\sigma) \to (v',v_*,\sigma)$ used in the cancellation 
lemma~\cite[Lemma~1]{ADVW:00} whose jacobian is $\cos^{-(N+\gamma)} (\theta/2)$ 
and is thus smaller than $2^{\frac{N+\gamma}{2}}$ for $\theta \in [0,\pi/2]$. 
For the second part we use the change of 
variable $(v,v_*,\sigma) \to (v,v_* ',\sigma)$ whose jacobian is 
$\sin^{-(N+\gamma)} (\theta/2)$, which is smaller than 
$2^{\frac{N+\gamma}{2}}$ for $\theta \in [\pi/2,\pi]$. Thus we get 
 \begin{multline*}
 \left|\int_{\R^N} Q^1 (g,f)(v) \phi(v) \, dv \right| \le \\ C_\Phi \, 
                    2^{\frac{N+\gamma}2} \, \|f\|_{W^{2,\infty}} \, 
                    \int_{\R^{2N} \times \ens{S}^{N-1}} 
                    dv \, dv_* \, d\sigma \, b(\cos \theta) \, 
                    \big( \sin (\theta/2) \big)^2 \, 
                    |v-v_*|^{2+\gamma} \, |g_*| \, |\phi| 
 \end{multline*}
and so, if $\Phi$ satisfy~\eqref{eq:Phi} with $2+\gamma \ge 0$ or 
if $\Phi$ satisfies~\eqref{eq:Phi:mol}, we get
 \begin{equation*}
 \left|\int_{\R^N} Q^1 (g,f)(v) \phi(v) \, dv \right| \le 
                    C_\Phi \, \frac{\|f\|_{W^{2,\infty}}}{2^{\frac{N+\gamma}2 +1}} \, 
                    \Bigg( \int_{\ens{S}^{N-1}} d\sigma \, 
                    b(\cos \theta) \, (1 - \cos \theta) \Bigg) 
                    \, \|g\|_{L^1 _{\tilde{\gamma}}} 
                    \, \| \phi \|_{L^1 _{\tilde{\gamma}}},  
 \end{equation*}
and if $\Phi$ satisfy~\eqref{eq:Phi} with $2+\gamma < 0$, then 
 \begin{multline*}
 \left|\int_{\R^N} Q^1 (g,f)(v) \phi(v) \, dv \right| \le \\ 
                    C_\Phi \, \frac{\|f\|_{W^{2,\infty}}}{2^{\frac{N+\gamma}2 +1}} \, 
                    \Bigg( \int_{\ens{S}^{N-1}} d\sigma \, 
                    b(\cos \theta) \, (1 - \cos \theta) \Bigg) 
                    \left[ \|g\|_{L^1 _{\tilde{\gamma}}} + \|g\|_{L^p} \right] 
                    \, \| \phi \|_{L^1 _{\tilde{\gamma}}}. 
 \end{multline*}
Since this holds for all $\phi \in L^1 _{\tilde{\gamma}}$, this yields the result by 
duality, with $\mbox{cst} = C_\Phi/(2^{\frac{N+\gamma}2 +1})$. 
\end{proof}

\section{Proof of the lower bound in the cutoff case}\label{sec:proof}
\setcounter{equation}{0}

In this section we shall prove Theorem~\ref{theo:main}. 
Since the collision operator is local in $t$ and $x$, 
the idea of the proof is to apply first Lemma~\ref{lem:up} and then Lemma~\ref{lem:spread} iterated 
on {\em each characteristic} of the free transport operator, in order first to obtain an upheaval 
point, and then to ``spread'' the minoration. It yields for each characteristic of the transport flow 
a Maxwellian lower bound with macroscopic velocity $\bar{v}$, temperature $\theta$ and 
density $\rho$ depending on the characteristic. 
Then a uniform control on $\bar{v}$, $\theta$ and $\rho$ yields the global Maxwellian lower bound.  
This control is based on the uniform bounds on the hydrodynamic quantities. 
The lower bound is also made uniform as $t \to +\infty$ thanks to these uniform bounds. 

The main tool is the Duhamel representation formula, written along the characteristics 
(which reduce to lines in the case of periodic boundary conditions):
 \begin{eqnarray}\nonumber
 && \forall \, t \in [0,T), \ \forall \, x \in \ens{T}^N, \ \forall \, v \, \in \R^N, \\ \label{eq:duham}
 && \quad f(t,x+vt, v) = f_0 (x,v) \exp\left(-\int_0 ^t L[f(s,x+vs,\cdot)](v) \, ds\right) \\ \nonumber
 && \qquad + \int_0 ^t \exp\left(- \int_s ^t L[f(s',x+vs',\cdot)](v) \, ds' \right) 
                 Q^+[f(s,x+vs,\cdot),f(s,x+vs,\cdot)](v) \, ds
 \end{eqnarray}
where $L$ was defined in~\eqref{eq:defL}.  
We define the concept of solution we shall use in the cutoff case, i.e. the concept of mild solutions 
(see~\cite[Section~5.3]{CIP:94} for an analogous definition).
 \begin{definition}\label{def:sol}
 Let $f_0$ be a measurable function non-negative almost everywhere on $\ens{T}^N \times \R^N$.  
 A measurable function $f=f(t,x,v)$ on $[0,T) \times \ens{T}^N \times \R^N$ is 
 a mild solution of the Boltzmann equation to the initial datum 
 $f_0 (x,v)$ if for almost every $(x,v)$ in $\ens{T}^N \times \R^N$:
  \[ t \mapsto L[f(t,x+vt,\cdot)](v), \ \ \ t \mapsto Q^+[f(t,x+vt,\cdot),f(t,x+vt,\cdot)](v) \]
 are in $L^1 _{\mbox{{\scriptsize {\em loc}}}} ([0;T))$, and for each $t \in [0,T)$, the 
 equation~\eqref{eq:duham} is satisfied and $f(t,x,v)$ is non-negative for almost every $(x,v)$.
 \end{definition}

 \begin{proposition}\label{prop:induc}
 Let $B = \Phi \,b$ be a collision kernel which satisfies~\eqref{eq:B}, with 
 $\Phi$ satisfying~\eqref{eq:Phi} or~\eqref{eq:Phi:mol}, and 
 $b$ satisfying~\eqref{eq:b} with $\nu < 0$.    
 Let $f(t,x,v)$ be a mild solution of the full Boltzmann equation in the torus  
 on some time interval $[0,T)$ ($T \in (0,+\infty]$), which satisfies 
  \begin{itemize}
  \item[(i)] assumption~\eqref{eq:hyp1} if $\Phi$ satisfies~\eqref{eq:Phi} with $\gamma \ge 0$ or 
  if $\Phi$ satisfies~\eqref{eq:Phi:mol},
  \item[(ii)] assumptions~\eqref{eq:hyp1} and~\eqref{eq:hyp2} if $\Phi$ satisfies~\eqref{eq:Phi} 
 with $\gamma \in (-N,0)$.  
 \end{itemize}
 Then for any fixed $\tau \in (0,T)$ and $x \in \ens{T}^N$, there exists some $R_0 >0$ 
 and some $\bar{v} \in B(0,R_0)$ such that
  \begin{equation}\label{eq:induc}
  \forall \, n \ge 0, \ \forall \, t \in \left[\tau - \frac{\tau}{2^{n+1}},\tau\right], \ \forall v \in \R^N, \ \ 
  f(t,x+vt,v) \ge a_n \, {\bf 1}_{B(\bar{v},\delta_n)}
  \end{equation}
 with the induction formulae
  \begin{equation*}\label{eq:reca}
  a_{n+1} = \mbox{{\em cst}} \; C_e \, \frac{a_n ^2 \delta_n ^{\gamma +N} \xi_n ^{N/2+1}}{2^{n+1}}
  \end{equation*}
  \begin{equation*}
  \delta_{n+1} = \sqrt{2} \, \delta_n (1 - \xi_n)
  \end{equation*}
 where $(\xi_n)_{n \ge 0}$ is any sequence in $(0,1)$, 
 $R_0 >0$, $a_0 >0$, $\delta_0 >0$, $C_e$  
 depend only on $\tau$, $B$, $\varrho_f$, $E_f$ and $H_f$ 
 (plus $L^{p_\gamma} _f$ if $\Phi$ satisfies~\eqref{eq:Phi} 
 with $\gamma \in (-N,0)$), and $\bar{v} \in B(0,R_0)$ depend on the same 
 quantities plus $x$. 
 \end{proposition}
\begin{proof}[Proof of Proposition~\ref{prop:induc}]
{\bf Step 1: Initialization}\\
We apply Lemma~\ref{lem:up} to the right-hand side member of the 
Duhamel representation iterated twice. More precisely, the equation~\eqref{eq:duham} 
yields on one hand
 \[ f(t,x+vt, v) \ge f_0 (x,v) \, e^{- C_L \, t \, \langle v \rangle^{\gamma^+}} \] 
and on the other hand 
 \[ f(t,x+vt, v) \ge \int_0 ^t \, e^{- C_L \, (t-s) \, \langle v \rangle^{\gamma^+}} 
                \, Q^+\big[f(s,x+vs,\cdot),f(s,x+vs,\cdot)\big](v) \, ds. \]
If we iterate the latter, we get 
 \begin{multline*}
 f(t,x+vt, v) \ge \int_0 ^t \, e^{- C_L \, (t-s) \, \langle v \rangle^{\gamma^+}} \\ 
 Q^+ \Big[ \Big( \int_0 ^s \, e^{- C_L \, (s-s') \, \langle v \rangle^{\gamma^+}} 
                \, Q^+ \big[ f(s',x+vs',\cdot),f(s',x+vs',\cdot) \big](\cdot) \, ds' \Big), 
                f(s,x+vs,\cdot) \Big](v) \, ds.
 \end{multline*}
Whenever $\varphi$ is some function on $\R^N$, we denote by $\varphi^{R_0}$ the 
truncation $\varphi \, {\bf 1}_{|v|\le R_0}$. We can bound from below by
 \begin{multline*}
 \forall \, v \in \R^N, \ |v| \le R_0 , \ \  
 f(t,x+vt, v) \ge Q^+ \Big[ Q^+ \big[ f_0 ^{R_0} (x,\cdot),f_0 ^{R_0} (x,\cdot) \big], 
                f_0 ^{R_0} (x,\cdot) \Big](v) \\ 
             \int_0 ^t e^{- C_L \, (t-s) \, R_0 ^{\gamma^+}} \, 
             e^{- C_L \, s \, R_0 ^{\gamma^+}} 
             \Big( \int_0 ^s  e^{- C_L \, (s-s') \, R_0 ^{\gamma^+}} \, 
              e^{- 2 \, C_L \, s' \, R_0 ^{\gamma^+}} \, ds' \Big) \, ds
 \end{multline*}             
and thus after some computation 
 \begin{multline*} 
 \forall \, v \in \R^N, \ |v| \le R_0 , \ \  
 f(t,x+vt, v) \ge  e^{- C_L \, t \, R_0 ^{\gamma^+}} 
                \, \frac{ \big( 1 - e^{- C_L \, t \, R_0 ^{\gamma^+}} \big)^2}
                {2 ( C_L \, R_0 ^{\gamma^+})^2} \\  
                Q^+ \Big[ Q^+ \big[ f_0 ^{R_0} (x,\cdot),f_0 ^{R_0} (x,\cdot) \big], 
                f_0 ^{R_0} (x,\cdot) \Big](v).  
 \end{multline*}
Then the use of Lemma~\ref{lem:up} concludes the initialization $n=0$ of the proof with 
 \[ a_0 = e^{- C_L \, \tau \, R_0 ^{\gamma^+}} 
                \, \frac{ \big( 1 - e^{- C_L \, (\tau/2) \, R_0 ^{\gamma^+}} \big)^2}
                {2 ( C_L \, R_0 ^{\gamma^+})^2} \, \eta_0\] 
where $\delta_0$, $R_0$, $\eta_0$ depend on $\varrho_f$, $E_f$, $H_f$ 
(as in the statement of Lemma~\ref{lem:up}) and $\bar{v}$ depend on the same 
quantities plus $x$ (via the function $f_0 ^{R_0} (x,\cdot)$).
\medskip \\
{\bf Step 2: Proof of the induction}\\
Now let us suppose that the induction property holds for $n$:
 \begin{equation*}
 \forall \, t \in \left[\tau - \frac{\tau}{2^{n+1}},\tau\right], \ \forall \, v \in \R^N, \ \ 
 f(t,x+vt,v) \ge a_n \, {\bf 1}_{B(\bar{v},\delta_n)}. 
 \end{equation*} 
The Duhamel representation yields the following lower bound
 \begin{multline*}
 \forall \, t \in \left[\tau - \frac{\tau}{2^{n+2}},\tau\right], \ \forall \, v \in \R^N, \ \ 
 f(t,x+vt,v) \ge  \\
 \int_{\tau-\frac{\tau}{2^{n+1}}} ^{\tau - \frac{\tau}{2^{n+2}}} 
 e^{- C_L \, (t-s) \, \langle v \rangle^{\gamma^+}}
 Q^+ \big( a_n \, {\bf 1}_{B(\bar{v},\delta_n)}, a_n \, {\bf 1}_{B(\bar{v},\delta_n)} \big) \, ds,  
 \end{multline*}
which easily leads to 
 \begin{multline*}
 \forall \, t \in \left[\tau - \frac{\tau}{2^{n+2}},\tau\right], \ \forall \, v \in \R^N, \ \ 
 f(t,x+vt,v) \ge \\ e^{- C_L \, \frac{\tau}{2^{n+1}} \, \langle v \rangle^{\gamma^+}} 
 \, \left( \frac{\tau}{2^{n+2}} \right) \, a_n ^2 \, 
 Q^+ \big( {\bf 1}_{B(\bar{v},\delta_n)}, {\bf 1}_{B(\bar{v},\delta_n)} \big). 
 \end{multline*}
Now the application of Lemma~\ref{lem:spread} gives 
 \begin{multline*}
 \forall \, t \in \left[\tau - \frac{\tau}{2^{n+2}},\tau\right], \ \forall \, v \in \R^N, \ \ 
 f(t,x+vt,v) \ge \\ \mbox{cst} \; e^{- C_L \, \frac{\tau}{2^{n+1}} \, \langle v \rangle^{\gamma^+}} 
 \, \left( \frac{\tau}{2^{n+2}} \right) \, a_n ^2 \, 
 \delta_n ^{N+\gamma} \, \xi_n ^{N/2-1} \, {\bf 1}_{B\big(\bar{v},\, \delta_n \sqrt{2} (1-\xi_n)\big)}
 \end{multline*}
and thus if $\delta_{n+1} = \delta_n \sqrt{2} (1-\xi_n)$, we get 
 \begin{multline*}
 \forall \, t \in \left[\tau - \frac{\tau}{2^{n+2}},\tau\right], \ \forall \, v \in \R^N, \ \ 
 f(t,x+vt,v) \ge \\ \mbox{cst} \; e^{- C_L \, \langle R_0 \rangle^{\gamma^+} \, 
 \frac{\tau}{2^{n+1}} \, \delta_{n+1} ^{\gamma^+}} 
 \, \left( \frac{\tau}{2^{n+2}} \right) \, a_n ^2 \, 
 \delta_n ^{N+\gamma} \, \xi_n ^{N/2-1} \, {\bf 1}_{B (\bar{v},\, \delta_{n+1} )}.
 \end{multline*}
As an easy induction shows, we have $\delta_n \le \delta_0 \, 2^{n/2}$ and so 
 \[ C_L \, \langle R_0 \rangle^{\gamma^+} \, 
 \frac{\tau}{2^{n+1}} \, \delta_{n+1} ^{\gamma^+} \le 
 C_L \; \langle R_0 \rangle^{\gamma^+} \, \frac{\tau}{2^{n+1}} \, (\delta_0 \, 2^{n/2})^{\gamma^+} \]
which is uniformly bounded from above since $\gamma \le 1$. Thus the exponential term 
 \[  e^{- C_L \, \langle R_0 \rangle^{\gamma^+} \, 
 \frac{\tau}{2^{n+1}} \, \delta_{n+1} ^{\gamma^+}} \]
is bounded from below uniformly by some constant $C_e > 0$. We deduce that 
 \begin{equation*}
 \forall \, t \in \left[\tau - \frac{\tau}{2^{n+2}},\tau\right], \ \forall \, v \in \R^N, \ \
 f(t,x+vt,v) \ge \mbox{cst} \; C_e \, \left( \frac{a_n ^2}{2^{n+1}} \right) \, a_n ^2 \, 
 \delta_n ^{N+\gamma} \, \xi_n ^{N/2-1} \, {\bf 1}_{B (\bar{v},\, \delta_{n+1} )}.
 \end{equation*}
This concludes the proof. 
\end{proof}

Now we can apply Proposition~\ref{prop:induc} along the characteristics in order to 
prove Theorem~\ref{theo:main}.
\begin{proof}[Proof of Theorem~\ref{theo:main}]
We shall divide the proof into three steps for the sake of clarity. Each step is embodied in a lemma. 
For these three lemmas we make the same assumptions on $B$ and $f$ as in Theorem~\ref{theo:main}.
\smallskip \\
{\bf Step 1: Choice of $(\xi_n)_{n \ge 0}$ and asymptotic behavior of $(a_n)_{n \ge 0}$}
 \begin{lemma}\label{lem:step1}
 For any $x \in \ens{T}^N$ and $\tau \in (0,T)$, there exists $\bar{v}(x) \in B(0,R_0)$ and 
 $\rho$, $\theta >0$ such that 
  \begin{equation} \label{eq:step1} 
  \forall \, v \in \R^N, \ \ \ 
  f(\tau,x+\tau v, v) \ge \rho \, \frac{e^{-\frac{|v-\bar{v}|^2}{2 \theta}}}{(2\pi \theta)^{N/2}}.
  \end{equation}
 The constants $R_0$, $\rho$, $\theta$ depend on $\tau$, 
 $\varrho_f$, $E_f$ and $H_f$ (and $L^{p_\gamma} _f$ if $\Phi$ satisfies~\eqref{eq:Phi} 
 with $\gamma <0$).
 \end{lemma}
\begin{proof}[Proof of Lemma~\ref{lem:step1}]
Let us now chose the sequence $(\xi_n)_{n \ge 0}$. The most natural choice is a geometrical 
sequence $\xi_n = \xi ^{n+1}$ for some $\xi \in (0,1)$. With this choice we can estimate 
the asymptotic behaviour of the sequence $(\delta_n)_{n \ge 0}$. Explicitely 
 \[ \delta_n = \delta_0 \, 2^{n/2} \, (1-\xi) \, (1 - \xi^2) \, 
               \cdots \, (1-\xi^n) = \delta_0 \, 2^{n/2} \, \Pi_{k=0} ^n (1-\xi^k) \] 
and thus as $\xi \in (0,1)$ one easily gets 
 \begin{equation*} 
 \delta_n \ge c_\delta \, 2^{n/2}
 \end{equation*}
where the constant $c_\delta$ depends on $\delta_0$ and $\xi$. 
It follows that 
 \begin{equation*}
 \forall \, n \ge 0, \ \forall \, t \in \left[\tau - \frac{\tau}{2^{n+1}},\tau\right], 
 \ \forall \, v \in B(\bar{v},c_\delta \, 2^{n/2}), \ \ 
 f(t,x+vt,v) \ge a_n.
 \end{equation*}
By plugging this into the expression of the Maxwellian distribution 
 \[  \rho \, \frac{e^{-\frac{|v-\bar{v}|^2}{2 \theta}}}{(2\pi \theta)^{N/2}} \] 
we deduce that a sufficient condition to obtain~\eqref{eq:step1} is the following 
lower bound on the coefficients $a_n$ appearing in the minoration~\eqref{eq:induc}: $a_n \ge \alpha ^{2^n}$   
for some $\alpha \in (0,1)$. Indeed the parameter $\theta$ can then be fixed such 
that 
 \[ e^{-\frac{c_\delta ^2}{2 \theta}} \le \alpha. \] 
Afterwards one can fix the parameter $\rho$ 
in order that 
 \[ a_0 \ge   \, \frac{\rho}{(2\pi \theta)^{N/2}} \] 
for $|v-\bar{v}| \le \delta_0$, which leads to~\eqref{eq:step1}.  

Let us prove this bound from below on 
the sequence $(a_n)_{n \ge 0}$. If one denotes 
 \[ \lambda_n = \frac{\delta_n ^{\gamma +N} \xi_n ^{N/2+1}}{2^{n+1}} \]
one gets explicitely 
 \begin{equation*}
 a_n = (\mbox{cst} \, C_e)^{2^n - 1} \, \left[ \lambda_{n-1} \, \lambda_{n-2} ^2 \, \cdots
 \, \lambda_{0} ^{2^{n-1}} \right] \, a_0 ^{2^n}. 
 \end{equation*}
As for the sequence $(\lambda_n)_{n \ge 0}$, we have $\lambda_n \ge \mbox{cst} \, \lambda ^{n}$ with 
 \[ \lambda = \frac{2^{(\gamma+N)/2} \, \xi^{N/2+1}}{2} \] 
and so  
 \[ a_n \ge (\mbox{cst} \, C_e)^{2^n-1} \, 
    \lambda^{\left[(n-1)2^0 + (n-2)2^1 + \cdots + 0 2^{n-1}\right]} \, a_0 ^{2^n}. \] 
If $\lambda >1$ the proof is clearly finished. If $\lambda \in (0,1)$ it remains to study 
the quantity 
 \[ A_n = \left[(n-1)2^0 + (n-2)2^1 + \cdots + 0 2^{n-1} \right] \]
An easy computation shows that $A_n = 2^n - (n+1)$ and so $A_n \le 2^n$. It yields $a_n \ge \alpha^{2^n}$ 
with $\alpha := \mbox{cst} \, C_e \, \lambda \, a_0$. 
\end{proof}
\noindent
{\bf Step 2: Uniformization of the spatial dependence} 
 \begin{lemma} \label{lem:step2} 
 For any $\tau \in (0,T)$, there exists $\rho', \theta' >0$ such that 
  \begin{equation} \label{eq:step2} 
  \forall \, x \in \ens{T}^N, \ \forall \, v \in \R^N, \ \ \ 
  f(\tau,x, v) \ge \rho' \, \frac{e^{-\frac{|v|^2}{2 \theta'}}}{(2\pi \theta')^{N/2}}.
  \end{equation}
 The constants $\rho'$, $\theta'$ depend on $\tau$, $\varrho_f$, 
 $E_f$ and $H_f$ (and $L^{p_\gamma} _f$ if $\Phi$ satisfies~\eqref{eq:Phi} 
 with $\gamma \in (-N,0)$).
 \end{lemma}
\begin{proof}[Proof of Lemma~\ref{lem:step2}] 
This step is straightforward: the right-member term in the estimate~\eqref{eq:step1} 
depends on the space variable $x$ only through $\bar{v} (x)$. However, 
as a consequence of Lemma~\ref{lem:up}, $\bar{v}$ is always included 
in the ball $B(0,R_0)$ for some radius $R_0$ depending only on the {\it a priori}  
bounds on the solution. Thus 
 \[ e^{-\frac{|v-\bar{v}|^2}{2 \theta}} \ge e^{-\frac{|v|^2}{\theta}} \, e^{-\frac{R_0 ^2}{\theta}} \]
and the proof is complete up to the choice of some new parameters $\rho', \theta'$: one can take 
$\theta' = \theta/2$ and  
 \[ \rho' = \rho \, \frac{e^{-\frac{R_0 ^2}{\theta}}}{2^{N/2}}. \]
\end{proof}

\noindent
{\bf Step 3: Uniformization of the time dependence}
 \begin{lemma} \label{lem:step3} 
 For any $\tau > 0$, there exists $\rho', \theta' >0$ such that 
  \begin{equation*} \label{eq:step3} 
  \forall \, t \in [\tau,T), \ \forall \, x \in \ens{T}^N, \ \forall \, v \in \R^N, \ \ \ 
  f(t,x, v) \ge \rho' \, \frac{e^{-\frac{|v|^2}{2 \theta'}}}{(2\pi \theta')^{N/2}}.
  \end{equation*}
 The constants $\rho'$, $\theta'$ depend on the {\em a priori} bounds on the solution. 
 \end{lemma}
\begin{proof}[Proof of Lemma~\ref{lem:step3}] 
Again this step is straightforward: one can check that the lower bound~\eqref{eq:step2} 
does not depend on the precise form of the solution $f(t,x,v)$ for $t \in [0,\tau]$, 
$x \in \ens{T}^N$, $v \in \R^N$, but only on the uniform bounds. It means 
that the same argument could be started not from $t=0$ anymore, but at any time 
(as long as the bounds used are uniform in time). As the lower bound appears after some time $\tau>0$ 
(arbitrary small), we get the lower bound for any time $t \ge \tau$ by making the proof 
start at $t-\tau$.
\end{proof}
This concludes the proof of Theorem~\ref{theo:main}. 
\end{proof}
 
\section{Proof of the lower bound in the non-cutoff case}\label{sec:proof:nc}
\setcounter{equation}{0}

In this section we shall prove Theorem~\ref{theo:main:nc}. 
Again we use the spreading property along each characteristic but 
we use the spreading property on the gain part of a truncated collision operator. 
The remaining part will be treated thanks to the $L^\infty$ 
estimates proved in Section~\ref{sec:tools}.  

We assume that $\nu \in [0,2)$ and we make the following splitting for any $\var \in (0,\pi/4)$:
 \[ Q  = Q^+ _\var - Q^- _\var + Q^1 _\var + Q^2 _\var \]
where $Q^+ _\var$ and $Q^- _\var$ are the usual Grad splitting for the collision 
operator with collision kernel 
 \[ B^S _\var := \Phi \, \left[ b \, {\bf 1}_{|\theta| \ge \var} \right] =: \Phi \, b^S _\var , \] 
and $Q^1 _\var$ and $Q^2 _\var$ are the splitting introduced in~\eqref{eq:split:nc} 
applied to the non-cutoff collision operator with collision kernel 
 \[ B^R _\var := \Phi \, \left[ b \, {\bf 1}_{|\theta| \le \var} \right] =: \Phi \, b^R _\var . \] 
For the sake of clarity the index $\var$ shall be recalled on each quantity that depends 
on this splitting. 

It is straightforward to check that $b^S _\var \ge \ell_b$ on $[\pi/4,3\pi/4]$, 
since $b^S _\var = b$ for $\theta \in [\pi/4,3\pi/4]$ and thus the constants given by the application 
of Lemma~\ref{lem:up} and Lemma~\ref{lem:spread} on $Q^+ _\var$ are uniform according to $\var$. 
Moreover we have
 \begin{equation}\label{eq:nm} 
 n_{b^S _\var} \sim_{\var \to 0} \frac{b_0}{\nu} \, \var^{-\nu}, \ \ \ \  
 m_{b^R _\var}  \sim_{\var \to 0} \frac{b_0}{2-\nu} \, \var^{2-\nu} 
 \end{equation}  
for $\nu \in (0,2)$ and 
 \begin{equation}\label{eq:nm:log} 
 n_{b^S _\var} \sim_{\var \to 0} b_0 \, |\log \var|, \ \ \ \  
 m_{b^R _\var}  \sim_{\var \to 0} \frac{b_0}{2-\nu} \, \var^2
 \end{equation}  
when $\nu =0$. 

The basic tool is the Duhamel formula written in the following way 
 \begin{eqnarray}\nonumber
 && \forall \, t \in [0,T), \ \forall \, x \in \ens{T}^N, \ \forall \, v \in \R^N, \\ \label{eq:duham:nc} 
 && \quad f(t,x+vt, v) = f_0 (x,v) \exp\left(-\int_0 ^t (S_\var +L_\var)[f(s,x+vs,\cdot)](v) \, ds\right) \\ \nonumber
 && \qquad   
      + \int_0 ^t \exp\left(- \int_s ^t (S_\var + L_\var)[f(s',x+vs',\cdot)](v) \, ds'\right) 
                \\ \nonumber
 && \qquad \qquad \qquad \qquad \qquad \qquad 
    \qquad (Q^+ _\var + Q^1 _\var) [f(s,x+vs,\cdot),f(s,x+vs,\cdot)](v) \, ds
 \end{eqnarray}
where $L_\var$ and $S_\var$ are the operators introduced in Section~\ref{sec:tools} 
corresponding respectively to $Q^- _\var$ and $Q^2 _\var$. 
We shall systematically use the $L^\infty$ estimates given by~\eqref{eq:LS} and 
Lemma~\ref{lem:Q1}, written in the following form
 \[ L_\var[f] \le C_f \, n_{b^S _\var} \, \langle v \rangle^{\gamma^+}, \ \ \ 
 S _\var [f] \le C_f \, m_{b^R _\var} \, \langle v \rangle^{\gamma^+}, 
 \ \ \ Q^1 _\var (f,f) \le C_f \, m_{b^R _\var} \, \langle v \rangle^{(2+\gamma)^+} \]
for a constant $C_f$ depending on the uniform bounds on $f$. 

Let us define the concept of mild solution we shall use in the non-cutoff case.
 \begin{definition}\label{def:sol:nc}
 Let $f_0$ be a measurable function non-negative almost everywhere on $\ens{T}^N \times \R^N$.  
 A measurable function $f=f(t,x,v)$ on $[0,T) \times \ens{T}^N \times \R^N$ is 
 a mild solution of the Boltzmann equation to the initial datum 
 $f_0 (x,v)$ if there exists $\var_0 > 0$ such that for all 
 $0 < \var < \var_0$, for almost every $(x,v)$ in $\ens{T}^N \times \R^N$:
  \[ t \mapsto Q^+ _\var [f(t,x+vt,\cdot),f(t,x+vt,\cdot)](v), \ \ 
     t \mapsto Q^1 _\var [f(t,x+vt,\cdot),f(t,x+vt,\cdot)](v), \] 
  \[ t \mapsto L_\var [f(t,x+vt,\cdot)](v),  \ \ t \mapsto S_\var [f(t,x+vt,\cdot)](v) \]
 are in $L^1 _{\mbox{\scriptsize {\em loc}}} ([0;T))$, and for each $t \in [0,T)$, the 
 equation~\eqref{eq:duham:nc} is satisfied and $f(t,x,v)$ is non-negative for almost 
 every $(x,v)$.
 \end{definition}

Let us prove the equivalent of Proposition~\ref{prop:induc} in the non cutoff case. 
Here we shall write the induction formula for a general sequence of time intervals $\Delta_n$. 

Indeed, on one hand at each step $n$ of the induction the spreading effect of the gain part $Q^+ _\var$ 
is now balanced by the perturbation $Q^1 _\var$, 
which imposes a careful choice of the splitting parameter $\var$ for each $n$ to get a 
lower bound on 
 \[ Q^+ _\var ( f,f) + Q^1 _\var ( f,f) \ge Q^+ _\var ( f,f)  - |Q^1 _\var (f,f)|. \]
This yields a sequence $(\var_n)_{n\ge0}$ going to $0$ as $n$ goes to infinity. 

On the other hand at each step $n$ of the induction 
we have the following action of $-Q^- _\var + Q^2 _\var$ along the characteristic in the estimate from below 
on the solution: 
 \[ e^{- C_f \, \big( m_{b^R _{\var_n}} + n_{b^S _{\var_n}} \big) \, 
 \big( \sum_{k \ge n+1} \Delta_k \big)
 \, \langle v \rangle^{\gamma^+}}. \]
It makes the lower bound decrease 
and this exponential term goes to $0$ when the splitting parameter $\var$ goes to $0$ 
(since $n_{b^S _{\var}}$ goes to infinity as $\var$ goes to $0$).  
That is why we shall choose time intervals $\Delta_n$ 
whose size decreases very fast to $0$ as $n$ goes to infinity, in order to limit the action 
of this part during a time interval. 

 \begin{proposition}\label{prop:induc:nc}
 Let $B = \Phi \,b$ be a collision kernel which satisfies~\eqref{eq:B}, with 
 $\Phi$ satisfying~\eqref{eq:Phi} or~\eqref{eq:Phi:mol}, and 
 $b$ satisfying~\eqref{eq:b} with $\nu \in [0,2)$.    
 Let $f(t,x,v)$ be a mild solution of the full Boltzmann equation in the torus  
 on some time interval $[0,T)$ ($T \in (0,+\infty]$), which satisfies 
  \begin{itemize}
  \item[(i)] assumptions~\eqref{eq:hyp1} and~\eqref{eq:hypnc} 
             if $\Phi$ satisfies~\eqref{eq:Phi} with $\gamma \ge 0$ or 
  if $\Phi$ satisfies~\eqref{eq:Phi:mol};
  \item[(ii)] assumptions~\eqref{eq:hyp1},~\eqref{eq:hyp2} and~\eqref{eq:hypnc} 
              if $\Phi$ satisfies~\eqref{eq:Phi} with $\gamma \in (-N,0)$.  
  \end{itemize}
 Then for any fixed $\tau \in (0,T)$ (small enough) and $x \in \ens{T}^N$, 
 any sequence $(\Delta_n)_{n \ge 0}$ of positive numbers such that $\sum_{n \ge 0} \Delta_n =1$, 
 there exists some $R_0 >0$ and $\bar{v} \in B(0,R_0)$ such that
  \begin{equation*}\label{eq:induc:nc}
  \forall \, n \ge 0, \ \forall \, t \in \Bigg[\Bigg(\sum_{k=0} ^n \Delta_k \Bigg) \tau,\tau \Bigg], 
  \ \forall v \in \R^N, \ \ 
  f(t,x+vt,v) \ge a_n \, {\bf 1}_{B(\bar{v},\delta_n)}.
  \end{equation*}
 The sequence $a_n$ satisfies the induction formula
  \begin{multline}\label{eq:reca:nc}
  a_{n+1} = \mbox{{\em cst}} \, \Delta_{n+1} \, 
  \exp \Bigg[ - \big[ C_f \, a_n ^2 \, \delta_n ^{N+\gamma-\tilde{\gamma}} 
                  \, \xi_n ^{(N/2 -1)} \big]^{-\frac{\nu}{2-\nu}} \, 
                  \Big( \sum_{k \ge n+1} \Delta_k \Big) \, \delta_{n+1} ^{\gamma^+} \Bigg] 
  \\ a_n ^2 \delta_n ^{\gamma +N} \xi_n ^{(N/2+1)}
  \end{multline}
 if $\nu \in (0,2)$ and 
  \begin{multline}\label{eq:reca:nc:log}
  a_{n+1} = \mbox{{\em cst}} \, \Delta_{n+1} \, 
  \exp \Bigg[ -\cst \; \log \big[ C_f \, a_n ^2 \, \delta_n ^{N+\gamma-\tilde{\gamma}} 
                  \, \xi_n ^{(N/2 -1)} \big] \, 
                  \Big( \sum_{k \ge n+1} \Delta_k \Big) \, \delta_{n+1} ^{\gamma^+} \Bigg] 
  \\ a_n ^2 \delta_n ^{\gamma +N} \xi_n ^{(N/2+1)}
  \end{multline}
 if $\nu =0$. The sequence $\delta_n$ satisfies the induction formula
  \begin{equation*}
  \delta_{n+1} = \sqrt{2} \, \delta_n (1 - \xi_n).
  \end{equation*}
 Here $(\xi_n)_{n \ge 0}$ is any sequence in $(0,1)$, 
 the constants $R_0>0$, $a_0>0$, $\delta_0>0$ and $C_f$ 
 depend only on $\tau$, $\varrho_f$, $E_f$, $E'_f$, $H_f$, $W_f$ 
 (plus $L^{p_\gamma} _f$ if $\Phi$ satisfies~\eqref{eq:Phi} 
 with $\gamma \in (-N,0)$), and $\bar v \in B(0,R_0)$ depends on the same quantities plus $x$. 
 \end{proposition}
\begin{proof}[Proof of Proposition~\ref{prop:induc:nc}]
In this proof we shall use estimates of Section~\ref{sec:tools} as well as several equations established 
in the proofs of Section~\ref{sec:proof}.
\medskip \\ 
{\bf Step 1: Initialization}\\
The initialization here is simpler than for the cutoff case since we assume 
some regularity on the solution, and thus we do not need the regularizing 
property of the iterated gain term. 
First we give a straightforward lemma:  
 \begin{lemma}\label{lem:up:nc}
 Let $g$ a non-negative function on $R^N$ such that $e_g$ and $w_g$ 
 are bounded, and $\rho_g$ satifies $0< \rho_g <+\infty$. 
 Then there are $R_0$, $\delta_0$, $\eta >0$ and $\bar{v} \in B(0,R_0)$
 such that 
  \[ g(v) \ge \eta \, {\bf 1}_{B(\bar{v},\delta_0)}, \] 
 where $R_0$, $\delta_0$, $\eta >0$ are explicit constants 
 depending on the upper bounds on $\rho_g$, $e_g$, $w_g$ and 
 the lower bound on $\rho_g$. 
 \end{lemma}
\begin{proof}[Proof of Lemma~\ref{lem:up:nc}]
Using the bound on the energy of $g$, the choice
 \[ R_0 = \sqrt{\frac{2 e_g}{\rho_g}} \] 
implies that 
 \[ \int_{|v|\le R_0} g(v) \, dv \ge \frac{\rho_g}{2}. \]
So there is $\bar{v} \in B(0,R_0)$ such that 
 \[ g(\bar{v}) \ge \frac{\rho_g}{2 \mbox{Vol}\big(B(0,R_0)\big)}. \]
As $w_g$ controls the Lipschitz norm, we have 
 \[ \forall \, v_1,v_2 \in \R^N, \ \ \ |g(v_1) - g(v_2)| \le w_g \, |v_1 -v_2| \] 
and thus if we take
 \[ \delta_0 =  \frac{\rho_g}{4  \mbox{Vol}\big(B(0,R_0)\big) w_g}, \ \ \ \ \ 
    \eta = \frac{\rho_g}{4 \mbox{Vol}\big(B(0,R_0)\big)}, \]
we get $g(v) \ge \eta \, {\bf 1}_{B(\bar{v},\delta_0)}$.  
\end{proof}
Now we fix $x \in \ens{T}^N$ and we deduce from the representation~\eqref{eq:duham:nc} that 
 \begin{multline*} 
 \forall \, t \in [0,\tau], \ \forall \, v \in \R^N, \ \ \ 
 f(t,x+vt, v) \ge f_0 (x,v) e^{-\int_0 ^t (S_\var +L_\var)[f(s,x+vs,\cdot)](v) \, ds} 
     \\ + \int_0 ^t e^{-\int_s ^t (S_\var + L_\var)[f(s',x+vs',\cdot)](v) \, ds'} 
                \, Q^1 _\var [f(s,x+vs,\cdot),f(s,x+vs,\cdot)](v) \, ds. 
 \end{multline*}
We apply Lemma~\ref{lem:up:nc} to the function $f_0(x,\cdot)$ to obtain
 \begin{multline*} 
 \forall \, t \in [0,\tau], \ \forall \, v \in \R^N, \ \ \ 
 f(t,x+vt, v) \ge \eta \, {\bf 1}_{B(\bar{v},\delta_0)} 
            e^{-\int_0 ^t (S_\var +L_\var)[f(s,x+vs,\cdot)](v) \, ds} 
     \\ + \int_0 ^t e^{-\int_s ^t (S_\var + L_\var)[f(s',x+vs',\cdot)](v) \, ds'} 
                \, Q^1 _\var [f(s,x+vs,\cdot),f(s,x+vs,\cdot)](v) \, ds. 
 \end{multline*} 
for some $\bar{v} \in B(0,R_0)$. Then we restrict the inequality on 
the ball $B(\bar{v},\delta_0) \subset B(0,R_0+\delta_0)$ and 
we use the estimates on $S_\var$, $L_\var$ and $Q^1 _\var$ to get 
 \begin{multline*} 
 \forall \, t \in [0,\tau], \ \forall \, v \in B(\bar{v},\delta_0), \ \ \  
 f(t,x+vt, v) \ge \eta \, {\bf 1}_{B(\bar{v},\delta_0)} 
            e^{- C_f \, \big( m_{b^R _\var} + n_{b^S _\var} \big) 
                     \, \langle v \rangle^{\gamma^+} \tau} 
     \\ - \tau \, C_f \, m_{b^R _\var} \, \langle v \rangle^{\tilde{\gamma}}.
 \end{multline*} 
Using the bounds on the velocity in the ball we get (up to modifying the constant 
$C_f$)
 \begin{multline*} 
 \forall \, t \in [0,\tau], \ \forall \, v \in B(\bar{v},\delta_0), \ \ \ 
 f(t,x+vt, v) \ge \eta \, {\bf 1}_{B(\bar{v},\delta_0)} 
            e^{- C_f \, \big( m_{b^R _\var} + n_{b^S _\var} \big) \, \tau} 
      -  \tau \, C_f \, m_{b^R _\var}. 
 \end{multline*}  
Then we assume (up to reducing $\tau$) that $\tau \le 1$, and we choose first $\var_0$ small 
enough such that 
 \[ C_f \, m_{b^R _{\var_0}} \le \frac{\eta}{4} \]
and then $\tau$ small enough such that 
 \[  e^{- C_f \, \big( m_{b^R _{\var_0}} + n_{b^S _{\var_0}} \big) \, \tau} \ge \frac12. \]
This shows that 
 \[  \forall \, t \in [0,\tau], \ \forall \, v \in B(\bar{v},\delta_0), \ \ \ 
 f(t,x+vt, v) \ge \frac{\eta}{4} \, {\bf 1}_{B(\bar{v},\delta_0)}, \]
which concludes the initialization with $\bar{v}$, $\delta_0$ and $\eta_0 = \eta/4$. 
\medskip \\
{\bf Step 2: Proof of the induction}\\
As for the proof of the induction, we proceed quite similarly as in the proof of Proposition~\ref{prop:induc}. 
We suppose that the $n$th step is satisfied: 
 \begin{equation*}
 \forall \, t \in \Bigg[\bigg(\sum_{k=0} ^n \Delta_k \Bigg) \tau, \tau \Bigg], 
 \ \forall \, v \in \R^N, \ \ 
 f(t,x+vt,v) \ge a_n \, {\bf 1}_{B(\bar{v},\delta_n)}. 
 \end{equation*}  
We use the following lower bound given by the Duhamel representation~\eqref{eq:duham:nc} and 
the estimate on $Q^1 _\var$:
 \begin{multline*}
 \forall \, t \in \Bigg[ \Bigg(\sum_{k=0} ^{n+1} \Delta_k \Bigg) \tau, \tau \Bigg], 
 \ \forall \, v \in \R^N, \\ 
 f(t,x+vt,v) \ge  \ 
 \int_{(\sum_{k=0} ^n \Delta_n) \tau} ^{t} 
 e^{- C_f \, \big( m_{b^R _\var} + n_{b^S _\var} \big) \, (t-s) \, \langle v \rangle^{\gamma^+}} \\
 \Bigg[ Q^+ _\var ( a_n \, {\bf 1}_{B(\bar{v},\delta_n)}, a_n \, {\bf 1}_{B(\bar{v},\delta_n)}) 
 -  \tau \, C_f \, m_{b^R _\var} \, \langle v \rangle^{\tilde{\gamma}} \Bigg] \, ds.
 \end{multline*}
Thus by applying Lemma~\ref{lem:spread} on $Q^+ _\var$, we obtain 
 \begin{multline*}
 \forall \, t \in \Bigg[\Bigg(\sum_{k=0} ^{n+1} \Delta_n \Bigg) \tau,\tau \Bigg], \ \forall \, v \in \R^N, \\ 
 f(t,x+vt,v) \ge \ 
 \int_{(\sum_{k=0} ^n \Delta_n) \tau} ^{t} 
 e^{- C_f \, \big( m_{b^R _{\var}} + n_{b^S _{\var}} \big) \, (t-s) \, \langle v \rangle^{\gamma^+}} \\
  \Bigg[ \cst \, a_n ^2 \, \delta_n ^{N+\gamma} \, \xi_n ^{(N/2 -1)} {\bf 1}_{B(\bar{v},\delta_{n+1})}
 -  \tau \, C_f \, m_{b^R _\var} \, \langle v \rangle^{\tilde{\gamma}} \Bigg]. 
 \end{multline*}
Then we restrict the inequality to the ball $B(\bar{v},\delta_{n+1})$ to obtain, using 
the bounds on the velocity and up to modifying the constant $C_f$:
 \begin{multline*}
 \forall \, t \in \Bigg[\Bigg(\sum_{k=0} ^{n+1} \Delta_n \Bigg) \tau,\tau \Bigg], 
 \ \forall \, v \in B(\bar{v},\delta_{n+1}), \\ 
 f(t,x+vt,v) \ge \ 
 \int_{(\sum_{k=0} ^n \Delta_n) \tau} ^{t} 
 e^{- C_f \, \big( m_{b^R _{\var}} + n_{b^S _{\var}} \big) \, (t-s) \, \delta_{n+1} ^{\gamma^+}} \\
  \Bigg[ \cst \, a_n ^2 \, \delta_n ^{N+\gamma} \, \xi_n ^{(N/2 -1)} {\bf 1}_{B(\bar{v},\delta_{n+1})}
 -  \tau \, C_f \, m_{b^R _\var} \, \delta_{n+1} ^{\tilde{\gamma}} \Bigg]. 
 \end{multline*}
Then (assuming $\tau \le 1$) we choose $\var = \var_n$ such that 
 \[ \tau \, C_f \, m_{b^R _\var} \, \delta_{n+1} ^{\tilde{\gamma}} 
    \le \frac12 \, \cst \, a_n ^2 \, \delta_n ^{N+\gamma} \, \xi_n ^{(N/2 -1)} \]
which is possible since $m_{b^R _\var} \to_{\var \to 0} 0$. More precisely by using 
the equivalent of $m_{b^R _\var}$ for $\var \sim 0$, simple computations 
show that we can take 
 \[ \var_n = \mbox{cst} \; 
    \Big[ C_f \, a_n ^2 \, \delta_n ^{N+\gamma-\tilde{\gamma}} \, \xi_n ^{(N/2 -1)} \Big]^{\frac{1}{2-\nu}} \]
where the constant $C_f$ is independent of $n$ and depends only on the uniform bounds on $f$. 

Then we restrict the time integration to the interval 
$[(\sum_{k=0} ^{n+1} \Delta_n) \tau, (\sum_{k=0} ^{n+1} \Delta_n) \tau]$ 
(since the integrand is non-negative) which yields 
 \begin{multline*}
 \forall \, t \in \Bigg[\Bigg(\sum_{k=0} ^{n+1} \Delta_n \Bigg) \tau,\tau \Bigg], 
 \ \forall \, v \in B(\bar{v},\delta_{n+1}), \\ 
 f(t,x+vt,v) \ge \cst \ 
 \int_{(\sum_{k=0} ^n \Delta_n) \tau} ^{(\sum_{k=0} ^{n+1} \Delta_n) \tau} 
 e^{- C_f \, \big( m_{b^R _{\var_n}} + n_{b^S _{\var_n}} \big) \, (t-s) \, \delta_{n+1} ^{\gamma^+}} 
 \, a_n ^2 \, \delta_n ^{N+\gamma} \, \xi_n ^{(N/2 -1)} {\bf 1}_{B(\bar{v},\delta_{n+1})} \\
 \ge \cst \ e^{- C_f \, \big( m_{b^R _{\var_n}} + n_{b^S _{\var_n}} \big) \, 
 \big( \sum_{k \ge n+1} \Delta_k \big) \, \delta_{n+1} ^{\gamma^+}} \, \Delta_{n+1} \, 
 \, a_n ^2 \, \delta_n ^{N+\gamma} \, \xi_n ^{(N/2 -1)} {\bf 1}_{B(\bar{v},\delta_{n+1})}.
 \end{multline*}

Finally the argument in the exponential is seen to be equivalent to 
 \[ \Big[ C_f \, a_n ^2 \, \delta_n ^{N+\gamma-\tilde{\gamma}} 
                \, \xi_n ^{(N/2 -1)} \Big]^{-\frac{\nu}{2-\nu}} \, 
                  \Bigg( \sum_{k \ge n+1} \Delta_k \Bigg) \, \delta_{n+1} ^{\gamma^+} \]
when $\nu \in (0,2)$ and 
 \[  - \cst \; \log \Big[ C_f \, a_n ^2 \, \delta_n ^{N+\gamma-\tilde{\gamma}} 
               \, \xi_n ^{(N/2 -1)} \Big] \, 
                  \Bigg( \sum_{k \ge n+1} \Delta_k \Bigg) \, \delta_{n+1} ^{\gamma^+} \]
when $\nu=0$ (for some new constant $C_f$ depending on the uniform bounds on $f$),  
which concludes the proof. 
\end{proof}
We are now able to conclude the proof of Theorem~\ref{theo:main:nc}.
\begin{proof}[Proof of Theorem~\ref{theo:main:nc}]
We only study the asymptotic behavior of the coefficients $a_n$. The 
two other steps of the proof (uniformization of the spatial and time dependences) 
are exactly similar as those in the proof of Theorem~\ref{theo:main}. 

We fix $\xi \in (0,1)$ and define $\xi_n = \xi^n$. We saw above that 
with this choice $\delta_n \sim \mbox{cst} \; 2^{n/2}$. 
First we deal with the case $\nu > 0$. 
Let us choose any $\kappa > 2+ 2\nu /(2-\nu)$, and take for the time intervals 
 \[ \Delta_{n+1} = \frac{\alpha^{\beta \kappa^n}}{\sum_{k \ge 0} \alpha^{\beta \kappa^{k-1}}} \] 
where $\alpha \in (0,1)$ and $2\nu /(2-\nu) < \beta < \kappa-2$. 
We shall establish by induction the lower bound
 \begin{equation}\label{eq:lba}
 a_n \ge \alpha^{\kappa^n}. 
 \end{equation}
One easily sees that this estimate~\eqref{eq:lba} 
implies that 
 \[ \forall \, v \in \R^N, \ \ \ 
  f(\tau,x +\tau v,v) \ge C_1 \, e^{-C_2 \, |v|^K} \]
for a suitable choice of $C_1, C_2 >0$ and
 \[ K = \frac{\log \kappa}{\log \sqrt{2}}, \]
and thus concludes the proof.

The initialization of the induction is made by choosing $\alpha$ such 
that $\alpha \le a_0$. Then we suppose the lower bound satisfied for $a_n$ and 
we show first that the argument of the exponential in~\eqref{eq:reca:nc} 
is uniformly bounded. A simple computation establishes that 
 \[ \sum_{k \ge n+1} \Delta_k \le \mbox{cst} \; \Delta_{n+1} = \mbox{cst} \; \alpha^{\beta \kappa^n}  \] 
where the constant is independent of $n$. Thus
 \begin{multline*}
 \Big[ a_n ^2 \, \delta_n ^{N+\gamma-\tilde{\gamma}} 
     \, \xi_n ^{(N/2 -1)} \Big]^{-\frac{\nu}{2-\nu}} \, 
    \Bigg( \sum_{k \ge n+1} \Delta_k \Bigg) \, \delta_{n+1} ^{\gamma^+} \\
    \le \left[ \frac{2^{\gamma^+ /2}}{\big( C_f 2^{(N+\gamma-\tilde{\gamma})/2)} 
           \xi^{(N/2 -1)} \big)^{\nu/(2-\nu)}} \right]^n \, 
                         \alpha^{\left(\beta -\frac{2 \nu}{2-\nu}\right) \kappa^n}
 \end{multline*}
and the right-hand side member of this inequality goes to $0$ as $n$ goes to infinity. 
So the exponential term is uniformly bounded from below by some constant $C_e >0$ depending 
on the uniform bounds on the solution $f$. 

So the induction formula~\eqref{eq:reca:nc} defining $a_n$ yields
 \[  a_{n+1} \ge \mbox{cst} \; C_e \, \Delta_{n+1} \, a_n ^2 \delta_n ^{\gamma +N} \xi_n ^{(N/2-1)} \]
and thus 
 \[ a_{n+1} \ge \mbox{cst} \; C_e \, \alpha^{(2+\beta) \kappa^n} \, 
                     \left[ 2^{(N+\gamma)/2)} \xi^{(N/2 -1)} \right]^n 
            \ge  \mbox{cst} \; \alpha^{\kappa^{n+1}}  \]
if $\alpha$ is small enough (using $\kappa > 2 + \beta$) and the induction is proved.

Now for the case $\nu =0$, we choose the following time intervals 
  \[ \Delta_{n+1} = \frac{\beta^n}{\sum_{k \ge 0} \beta^{k-1}} \] 
where $\beta \in (0,1)$. We shall establish by induction the lower bound
 \begin{equation*}
 a_n \ge \alpha^{2^n}
 \end{equation*}
which implies (as in the proof of Theorem~\ref{theo:main}) that 
 \[ \forall \, v \in \R^N, \ \ \ 
  f(\tau,x +\tau v,v) \ge C_1 \, e^{-C_2 \, |v|^2} \]
for a suitable choice of $C_1, C_2 >0$ and thus concludes the proof.

We suppose the lower bound satisfied for $a_n$ and we show first that 
the argument of the exponential in~\eqref{eq:reca:nc:log} is uniformly bounded. We have 
 \[ \sum_{k \ge n+1} \Delta_k \le \mbox{cst} \; \Delta_{n+1} = \mbox{cst} \; \beta^n  \] 
where the constant is independent of $n$, and so 
 \begin{multline*}
 \left| \log \Big[ C_f \, a_n ^2 \, \delta_n ^{N+\gamma-\tilde{\gamma}} \, \xi_n ^{(N/2 -1)} \Big] \, 
                  \Bigg( \sum_{k \ge n+1} \Delta_k \Bigg) \, \delta_{n+1} ^{\gamma^+} \right| \\
     \le |\log C_f | + \mbox{cst} \, 2^{n(\gamma^+/2+1)} \, \beta^n + 
     \mbox{cst} \, n \, 2^{n\gamma^+/2} \, \beta^n
 \end{multline*}
which goes to $0$ if $\beta$ is taken small enough. So the exponential term 
is uniformly bounded from below by some constant $C_e > 0$ depending on the 
uniform bounds on the solution $f$. 

The induction formula~\eqref{eq:reca:nc:log} defining $a_n$ yields 
  \[  a_{n+1} \ge \mbox{cst} \; C_e \, \Delta_{n+1} \, a_n ^2 \delta_n ^{\gamma +N} \xi_n ^{(N/2-1)} \]
and thus 
 \[ a_{n+1} \ge \mbox{cst} \; C_e \, a_n ^2 \, 
                     \left[ 2^{(N+\gamma)/2)} \xi^{(N/2 -1)} \beta \right]^n. \]
Then if we denote 
 \[ \lambda = 2^{(N+\gamma)/2)} \xi^{(N/2 -1)} \beta \] 
a similar computation as in the proof of Proposition~\ref{prop:induc} gives 
 \[ a_n \ge 
  (\mbox{cst} \, C_e)^{2^n-1} \, \lambda^{\left[(n-1)2^0 + (n-2)2^1 + \cdots + 0 2^{n-1}\right]} \, a_0 ^{2^n} 
     \ge \left( \mbox{cst} \, C_e \, \lambda \, a_0 \right)^{2^n} \]
and thus $a_n \ge \alpha^{2^n}$ if one takes $\alpha \le \mbox{cst} \, C_e \, \lambda \, a_0$ and the 
induction is proved. This concludes the proof. 
\end{proof}

\section{Application to the existing Cauchy theories}\label{sec:app}
\setcounter{equation}{0}

In this section we shall use Theorem~\ref{theo:main} and Theorem~\ref{theo:main:nc} 
to study solutions which have been constructed by previous authors, by connecting these 
theorems to some existing results in the Cauchy theory of the Boltzmann equation.

First we give a theorem which summarizes the situation in the 
spatially homogeneous setting for cutoff potentials 
in cases where the collision kernel does not present a singularity 
for vanishing relative velocity (a case which is not so well understood 
and for which $L^p$ estimates have not yet been derived). 

 \begin{theorem}\label{theo:app:homct}
 Let $B =\Phi \, b$ be a collision kernel satisfying assumptions~\eqref{eq:B}, with 
 $\Phi$ satisfying assumption~\eqref{eq:Phi} with $\gamma \ge 0$ or~\eqref{eq:Phi:mol}, 
 and $b$ satisfying~\eqref{eq:b} with $\nu < 0$. 
 Let $f_0$ be a nonnegative initial condition on $\R^N _v$ with finite mass, energy.  
 Then 
  \begin{itemize}
  \item[(i)] there exists a unique solution $f(t,v)$ with constant mass and energy to the spatially 
  homogeneous Boltzmann equation, defined for all times;
  \item[(ii)] if $f_0$ has finite entropy, then the entropy of the solution remains uniformly 
  bounded and the solution satisfies
   \begin{equation*}
   \forall \, t > 0, \ \forall \, v \in \R^N, \ \ \ 
   f(t,v) \ge \rho(t) \, \frac{e^{-\frac{|v|^2}{2 \theta(t)}}}{(2\pi \theta(t) )^{N/2}}.
   \end{equation*}
  The constants $\rho(t)$, $\theta(t) >0$ are explicit and depend on the mass, energy and 
  entropy of $f_0$; they are uniform for $t \to +\infty$ but not necessarily for $t \to 0$. 
  \end{itemize}
 \end{theorem}

\Remark  Let us sketch briefly how 
it is possible to relax the assumption of the 
boundedness of the entropy of the initial datum in point (ii)  
in the case $\gamma >0$ in dimension $3$. Indeed 
Mischler and Wennberg~\cite[Lemma~2.1]{MiscWenn:exun:99} proved in this case that  
 \[ g = Q^+ (Q^+(f,f),f) \] 
is uniformly integrable, with constants depending on the $L^1 _2$ norm of $f$. 
The bound on the entropy is only used in the obtaining of the upheaval point in 
Lemma~\ref{lem:up}, whose proof requires the uniform integrability of the function. 
But in the initialization step of Proposition~\ref{prop:induc}, 
it is possible, up some tricky computations, to obtain by iterating twice more 
the Duhamel representation
 \begin{multline*} 
 \forall \, v \in \R^N, \ |v| \le R_0 , \ \  
 f(t,x+vt, v) \ge C_{T,R_0,B} \, 
                Q^+ \Big[ Q^+ \big[ g_0 ^{R_0} (x,\cdot),g_0 ^{R_0} (x,\cdot) \big], 
                g_0 ^{R_0} (x,\cdot) \Big](v)
 \end{multline*}
where $g_0 = Q^+ (Q^+(f_0,f_0),f_0)$. As $g_0$ is uniformly integrable (with explicit bounds) 
by the result of Mischler and Wennberg above, one can apply Lemma~\ref{lem:up} to 
 \[ Q^+ \Big[ Q^+ \big[ g_0 ^{R_0} (x,\cdot),g_0 ^{R_0} (x,\cdot) \big], 
                g_0 ^{R_0} (x,\cdot) \Big](v) \]
to get the upheaval point and the rest of the proof is unchanged.
\medskip
 
\begin{proof}[Proof of Theorem~\ref{theo:app:homct}]
Let us prove (i): In the case  $\Phi$ satisfies assumption~\eqref{eq:Phi} with $\gamma > 0$, 
the existence and uniqueness in $L^1 _2$ (for solutions with 
non-increasing energy) are proved in~\cite{MiscWenn:exun:99}. In the case $\gamma=0$ or 
$\Phi$ satisfies assumption~\eqref{eq:Phi:mol} with $\gamma \in (-N,0)$, 
existence and uniqueness in $L^1 _2$ can be deduced from Arkeryd~\cite{Arke:I+II:72}: in this case 
the collision operator is a bounded bilinear operator in $L^1$, which implies the uniqueness, and 
the global existence is proved by the monotonicity argument from~\cite{Arke:I+II:72}. 
For (ii): In all these cases the mass and energy are conserved. By the $H$ theorem, if the 
entropy of the initial datum is bounded, then it remains bounded uniformly for all times. 
Thus the solution satisfies~\eqref{eq:hyp1} and one can apply Theorem~\ref{theo:main} 
and concludes the proof. 
\end{proof}
  
Now we give a theorem for non-cutoff mollified hard potentials collision kernels, 
using a recent result of Desvillettes and Wennberg~\cite{DesvWenn:nc:pp}. 
Here ${\cal S}(\R^N _v)$ denotes the Schwartz space of the functions with all derivatives 
bounded and decreasing faster at infinity than any inverse of polynomial. 
 \begin{theorem}\label{theo:app:homnc}
 Let $B =\Phi \, b$ be a collision kernel satisfying assumptions~\eqref{eq:B}, with 
 $\Phi$ satisfying assumption~\eqref{eq:Phi:mol} with $\gamma > 0$ and $C^\infty$, 
 and $b$ satisfying~\eqref{eq:b} with $\nu \in (0,2)$. 
 Let $f_0$ be a nonnegative initial condition on $\R^N _v$ with finite mass, energy and entropy.  
 Then 
  \begin{itemize}
  \item[(i)] there exists a solution $f$ to the spatially homogeneous Boltzmann equation with 
  constant mass and energy and uniformly bounded entropy, 
  defined for all times and belonging to $L^\infty ([t_0,+\infty),{\cal S}(\R^N _v))$ for any 
  $t_0 > 0$; 
  \item[(ii)] this solution satisfies 
   \begin{equation*}
   \forall \, t > 0, \ 
   \forall \, v \in \R^N, \ \ \ 
   f(t,v) \ge C_1 (t) \, e^{-C_2 (t) \, |v|^K}
   \end{equation*}
  for any exponent $K$ such that 
   \[ K > 2 \, \frac{\log \left(2 + \frac{2 \nu}{2 - \nu}\right)}{\log 2}. \] 
  The constants $C_1 (t)$, $C_2 (t) >0$ are explicit and depend on the mass, energy, entropy of 
  $f_0$ and $K$; they are uniform for $t \to +\infty$ but not necessarily for $t \to 0$. 
  \end{itemize}
 \end{theorem}

\begin{proof}[Proof of Theorem~\ref{theo:app:homnc}]
First let us prove (i): Our assumptions on $B$ imply the assumptions~\cite[Assumption~2]{DesvWenn:nc:pp} 
on the collision kernel $B$, namely $B=\Phi \, b$ with  
$\Phi$ a smooth and strictly positive function such that 
$\Phi(z) \sim_{|z| \to +\infty} |z|^\gamma$ with $\gamma \in (0,1]$, and $b$ 
such that $b(\cos \theta) \sim_{\theta \to 0} \mbox{cst} \; \theta^{-(N-1)-\nu}$ with 
$\nu >0$. Concerning the initial datum our assumptions are exactly those 
of~\cite[Assumption~1]{DesvWenn:nc:pp}. 
So we can apply~\cite[Theorem~1]{DesvWenn:nc:pp} to prove the existence of a solution, 
lying in $L^\infty ([t_0,+\infty),{\cal S}(\R^N _v))$ for any $t_0 > 0$.  
For (ii): 
The explicit bound $L^\infty ([t_0,+\infty),{\cal S}(\R^N _v))$ for any $t_0 >0$ 
immediately implies the uniform bounds~\eqref{eq:hyp1} and~\eqref{eq:hypnc}. 
Thus one can apply Theorem~\ref{theo:main:nc} to obtain the lower 
bound for $t \ge t_0+\tau$. As $t_0$ and $\tau$ are arbitrarily small 
this concludes the proof.
\end{proof}

For spatially inhomogeneous solutions we can apply 
our results to the solutions near the equilibrium in a torus constructed 
by Ukai (see~\cite{Ukai:FBE:74,Ukai:FBE:76} and~\cite[Section~7.6]{CIP:94}) 
for hard spheres and Guo~\cite{Guo:mous:01} for soft potentials. 
For the sake of clarity, we do not explicit in full details the functional settings in which these 
solutions are constructed and refer to the above-mentioned references for more precise 
definitions. 

 \begin{theorem} \label{theo:app:inhom}
 Let $B =\Phi \, b$ be a collision kernel satisfying assumptions~\eqref{eq:B}, with 
  \begin{itemize}
  \item[a-] $\Phi$ satisfying assumption~\eqref{eq:Phi} with 
  $\gamma = 1$ and $b=1$ (solutions of Ukai) or 
  \item[b-] $\Phi$ satisfying assumption~\eqref{eq:Phi} with 
  $\gamma < 0$ and $b$ satisfying~\eqref{eq:b} with $\nu < 0$ (solutions of Guo).
 \end{itemize}
 Let $f_0 = M + M^{1/2} h_0$ ($M$ is the global Maxwellian equilibrium) 
 be a nonnegative initial condition on $\ens{T}^N _x \times \R^N _v$ such that 
  \[ \left\|h_0\right\|^2 _{H^{s,q}} := \sum_{|i|+|j|\le s} \
             \left\|h_0 \langle v \rangle^q \right\|^2 _{L^2(\ens{T}^N \times \R^N)} 
                       \le \epsilon_0 \]
 with $s, q, \epsilon_0 >0$. Then 
  \begin{itemize}
  \item[(i)] for any $s',q'>0$, if $s,q$ are large enough and $\epsilon_0$ is small enough, 
  there exists a unique solution $f$ to the full Boltzmann equation in $H^{s',q'}$, 
  defined for all times, with uniform bound in $H^{s',q'}$ depending on $\epsilon_0$; 
  \item[(ii)] this solution satisfies
   \begin{equation*}
   \forall \, t > 0, \ 
   \forall \, x \in \ens{T}^N, \ 
   \forall \, v \in \R^N, \ \ \ 
   f(t,x,v) \ge \rho(t) \, \frac{e^{-\frac{|v|^2}{2 \theta(t)}}}{(2\pi \theta(t) )^{N/2}}.\
   \end{equation*}
  The constants $\rho(t)$, $\theta(t) >0$ depend on $\epsilon_0$; 
  they are uniform for $t \to +\infty$ but not necessarily for $t \to 0$. 
  \end{itemize}
 \end{theorem}
 
\Remarks 1. Most probably the solutions of Ukai could extend to any cutoff hard potentials, even 
if they were constructed for hard spheres. Anyhow the method of Guo would probably recover 
the result for any cutoff hard potentials as well. 
\smallskip

2. The proof of Ukai uses the spectral gap of the linearized collision operator for 
hard spheres, which was known to exist since Grad~\cite{Grad:58}. This spectral 
gap was obtained by non-constructive method (essentially the Weyl's Theorem about 
compact perturbation of the essential spectrum). However explicit estimates 
on this spectral gap were recently obtained in~\cite{BaMo:sg:03}, which is a step 
forward into a constructive theory. Nevertheless at now the proofs of Ukai 
(and Guo for soft potentials) still do not provide constructive bounds.
\smallskip

3. One can apply the result of convergence to equilibrium of Desvillettes and 
Villani~\cite{DV:eqEB:03} to the solutions of Theorem~\ref{theo:app:inhom}, which 
do satisfy every assumption of~\cite[Theorem~2]{DV:eqEB:03}. Thus they converges 
almost exponentially to equilibrium, i.e 
 \[ \|f_t - M\|_{L^1 _{x,v}} \le C_\alpha \, t^{-1/\alpha} \]
where $\alpha$ can be taken as big as wanted when $s$ and $q$ are large enough, 
and the constant $C_\alpha$ is explicit according to 
the uniform regularity bounds on $f$ and the constants in the lower bound. 
For the solutions of Ukai in the hard spheres case, 
this result is weaker than the exponential convergence 
to equilibrium already proved by Ukai, but the convergence 
to equilibrium was unknown for the solutions of Guo in the case of soft potentials. 
\smallskip

4. Let us explain how to skip the ``upheaval step'' of the proof for solutions 
near the equilibrium. Indeed in the case of perturbative solutions, 
which are $L^\infty$ close to a global Maxwellian distribution, 
Lemma~\ref{lem:up} can be bypassed by a more direct argument: 
up to reduce the neightborhood of the Maxwellian distribution for the existence theory 
(i.e. reducing $\epsilon_0$), the uniform $L^\infty$ control of smallness on $h$ yields 
 \[ \forall \, t \in \R_+, \ \forall \, x \in \ens{T}^N, \ \forall \, v \in \R^N, \ \ \ 
    f(t,x,v) \ge \eta_0 \, {\bf 1}_{B(0,\delta_0)} \]
for some $\eta_0 >0$ and $\delta_0 >0$. 
\medskip

\begin{proof}[Proof of Theorem~\ref{theo:app:inhom}]
One can check easily that for $s$ and $q$ large enough and $\epsilon_0$ small enough 
in the assumptions on the initial data, $f_0$ satisfies the assumptions 
of~\cite[Theorem~7.6.2]{CIP:94} in the case of Ukai's solutions, and the assumptions 
of~\cite[Theorem~1]{Guo:mous:01} for Guo's solutions. Then for $s$ large enough, 
the uniform smoothness estimates on the solution imply in both cases a bound on $\|h_t\|_{L^\infty _{x,v}}$, 
uniform for all times, and going to $0$ as $\epsilon_0 \to 0$. 
Thus if one take $\epsilon_0$ small enough such that 
 \[  \left( \sup_{t \ge 0} \|h_t\|_{L^\infty _{x,v}} \right) \, \int_{\R^N _v} \sqrt{M(v)} \, dv 
      \le \frac12 \,  \int_{\R^N _v} M(v) \, dv \]
one immediately get 
 \[ \varrho_f \ge \frac12 \,  \int_{\R^N _v} M(v) \, dv  > 0 \]
for all times. The uniform upper bounds on the local energy and local entropy 
(and local $L^p$ bound) follow 
from the uniform regularity bounds on the solution. Thus the solution 
satisfies~\eqref{eq:hyp1} and~\eqref{eq:hyp2} and one can apply Theorem~\ref{theo:main} 
and conclude the proof. 
\end{proof}

Finally, let us say a few words about other Cauchy theories. 

In the spatially inhomogeneous setting, one could apply Theorem~\ref{theo:main} to 
solutions for small time constructed in~\cite{KaSh:EB:78} (in the cutoff case). 
We did not detail this application since we were more interested with global solutions. 

One could also apply Theorem~\ref{theo:main} to 
the global weakly inhomogeneous solutions (for cutoff hard potentials) constructed by Arkeryd, 
Esposito and Pulvirenti in~\cite{AEP:EB:87}, which would give a Maxwellian 
lower bound on these solutions (uniform as $t \to +\infty$). Note that the uniform bounds 
on the solution obtained in~\cite{AEP:EB:87} do not seem to be constructive. 
As a consequence the Maxwellian lower bound given by Theorem~\ref{theo:main} 
would not have constructive constants.   

Concerning the global solutions in the whole space $\R^N$ near the vacuum 
constructed by Kaniel, Illner and Shinbrot 
(cf.~\cite{IlSh:EB:84,MiPe:inhom:97,Goud:maxw:97} 
or~\cite[Section~5.2]{CIP:94}), 
a lower bound on the solution $f(t,x,v)$ cannot be uniform in space since $f(t,\cdot,\cdot)$ is 
integrable on $\R^N _x \times \R^N _v$, and it cannot be uniform as $t$ goes to infinity 
since the solution goes to $0$ as $t$ goes to infinity for every $(x,v)$ such that 
$v \neq 0$ (see~\cite[Theorem~5.2.2]{CIP:94}). 
Our method could not apply, and is more adapted to evolution problems in bounded domains. 
We note that for these solutions in the whole space, 
in some cases a bound from below on the solutions by a ``travelling Maxwellian'' 
 \[ \forall \, t \ge 0, \ \forall \, x \in \R^N, \ \forall \, v \in \R^N, \ \ \ 
    f(t,x,v) \ge C(t) \, e^{-\beta \, |x-tv|^2} \, e^{-\alpha \, |v|^2} \] 
(where $\alpha > 0$ and $\beta > 0$ are absolute constants, and $C(t) > 0$ is a constant depending 
on time) can replace our method to provide a lower bound (see for instance~\cite{Goud:maxw:97}, 
and also, in the same spirit~\cite{Lu:99}). 

\medskip
\noindent
{\bf{Acknowledgment}}: Support by the European network HYKE, 
funded by the EC as contract HPRN-CT-2002-00282, is acknowledged. 
\medskip


\begin{flushright} \signcm \end{flushright}

\end{document}